\documentclass[12pt,thmsa]{article}%
\usepackage{amsmath}
\usepackage{amssymb}
\usepackage{amsfonts}
\usepackage{graphicx}%
\newcommand{\slg}{\mbox{\bfseries\slshape g}}
\begin{document}

\title{Geometric Algebras and Extensors}
\author{{\footnotesize V. V. Fern\'{a}ndez}$^{1}${\footnotesize ,} {\footnotesize A.
M. Moya}$^{2}$ {\footnotesize and W. A. Rodrigues Jr}.$^{1}${\footnotesize \ }\\$^{1}\hspace{-0.1cm}${\footnotesize Institute of Mathematics, Statistics and
Scientific Computation}\\{\footnotesize \ IMECC-UNICAMP CP 6065}\\{\footnotesize \ 13083-859 Campinas, SP, Brazil }\\{\footnotesize e-mail: walrod@ime.unicamp.br \ virginvelfe@accessplus.com.ar}\\{\footnotesize \ }$^{2}${\footnotesize Department of Mathematics, University
of Antofagasta, Antofagasta, Chile} \\{\footnotesize e-mail: mmoya@uantof.cl}}
\maketitle

\begin{abstract}
This is the first paper in a series (of four) designed to show how to use
geometric algebras of multivectors and extensors to a novel presentation of
some topics of differential geometry which are important for a deeper
understanding of geometrical theories of the gravitational field. In this
first paper we introduce the key algebraic tools for the development of our
program, namely the euclidean geometrical algebra of multivectors
$\mathcal{C}\ell(V,G_{E})$ and the theory of its deformations leading to
metric geometric algebras $\mathcal{C}\ell(V,G)$ and some special types of
extensors. Those tools permit obtaining, the remarkable golden formula
relating calculations in $\mathcal{C}\ell(V,G)$ with easier ones in
$\mathcal{C}\ell(V,G_{E})$ ( e.g., a noticeable relation between the Hodge
star operators associated to $G$ and $G_{E}$). Several useful examples are
worked in details fo the purpose of transmitting the \textquotedblleft tricks
of the trade\textquotedblright.

\end{abstract}
\tableofcontents

\section{Introduction}

This is the first paper in series of four, designed to show how Clifford
(geometric) algebra methods can be conveniently used in the study of
differential geometry and geometrical theories of the gravitational field. It
dispenses the use of fiber bundle theory and is indeed a very powerful and
economic tool for performing sophisticated calculations. This first paper deal
with the algebraic aspects of the theory, namely Clifford algebras and the
theory of extensors. Our presentation is self contained and serves besides the
purpose of fixing conventions also the one of introducing a series of
\textquotedblleft tricks of the trade\textquotedblright\ (not found easily
elsewhere) necessary for quickly and efficient computations. The other three
papers develop in a systematic way a theory of \textit{multivector} and
\textit{extensor} fields \cite{3} and their use in the differential geometry
of manifolds of arbitrary topology \cite{4,7}. There are many novelties in our
presentation of differential geometry, in particular the way we introduce the
concept of deformation of geometric structures, which is discussed in detail
in \cite{4} and which permit us \cite{4,7} to relate some well distinct
geometric structures on a given manifold. Moreover, the method permit also to
solve problems in one given geometry in terms of an eventually simple one, and
here is a place where our theory may be a useful one for the study of
geometrical theories of the gravitational field.

The main issues discussed in the present paper are the constructions of the
Euclidean and the metric \textit{geometric} (or Clifford) algebras of
\emph{multivectors}, denoted $\mathcal{C\ell}(V,G_{E})$ and $\mathcal{C\ell
}(V,G)$ which can be associated to a real vector space $V$ of dimension $n$
once we equip $V$ \ respectively with an Euclidean $(G_{E})$ and an arbitrary
$(G)$ non degenerated metric of signature $(p,q)$ with $p+q=n$. The Euclidean
geometrical algebra is the \textit{key} tool for performing almost all
calculations of the following papers. It plays in our theory a role
\textit{analogous} to matrix algebra in standard presentations of linear
algebra. Metric geometric algebras are introduced as \emph{deformations} of
the Euclidean geometric algebra. This is conveniently explored with the
introduction of the concept of a \textit{deformation}\emph{\ extensor}
associated to a given metric extensor and and by proving the remarkable
\emph{golden formula}. \textit{Extensors} are a new kind of geometrical
objects which play a crucial role in the theory presented in this series and
in what follows the basics of their theory is described. These objects have
been apparently introduced by Hestenes and Sobczyk in \cite{8} and some
applications of the concept appears in \cite{9}, but a rigorous theory was
lacking until recently\footnote{More details on the theory of extensors may be
found in \cite{10}.}. It is important to observe that in \cite{8} the
preliminaries of the geometric calculus have been applied to the study of the
differential geometry of \emph{vector} manifolds. However, as admitted in
\cite{11} there are some problems with such approach. In contrast, our
formulation applies to manifolds of an arbitrary topology and is free of the
problems that paved the construction in \cite{8}. On the following papers
dealing with differential geometry, the concept of a \textit{canonical vector
space} associated to an open set of an arbitrary manifold is introduced and
the main ingredients of differential geometry, like connections and their
torsion and curvature extensors are introduced using the Euclidean geometric
calculus in the canonical vector space. After that metric extensors are
introduced and the concept of \textit{deformed geometries} relative to a given
geometry is introduced. An \emph{intrinsic} Cartan calculus is developed
together with some other topics that are ready to be used in geometrical
theories of the gravitational field, as, we hope, the reader will convince
himself consulting the other papers in the series.

As to the explicit contents of the present paper they are summarized as
follows. Section 2 is dedicated to the construction of the Euclidean and
metric geometric algebras. The concepts of \ \textit{multivectors}
(homogeneous and non homogeneous ones) and their \textit{exterior algebra }is
briefly recalled\textit{. }Next we introduce the \textit{scalar product} of
multivectors and the concepts of right and left \textit{contractions} and
\textit{interior} algebras, and then give a definition of a general
\emph{real} Clifford (or geometrical) algebra $\mathcal{C\ell}(V,G)$ of
multivectors associated to a pair $(V,G)$. We fix an Euclidean metric $G_{E}$
in $V$ and use $\mathcal{C\ell}(V,G_{E})$, the Euclidean geometric algebra as
our basic tool for performing calculations, and study some relations between
$\mathcal{C\ell}(V,G)$ and $\mathcal{C\ell}(V,G_{E})$ \footnote{As additional
references to several aspects of the theory of the theory of Clifford
algebras, that eventually may help the interested reader, we quote
\cite{1,10,6,12,13,14,rod041,rodoliv2006}}. Section 3 is dedicated to an
introduction to the theory of extensors, with emphasis on the properties of
some special kind of extensors that will be used in our approach to
differential geometry in the next papers of the series. Of special interest is
the golden formula which permit us to make calculations related to
$\mathcal{C\ell}(V,G)$ using the simple algebra $\mathcal{C\ell}(V,G_{E})$
once we determine a gauge extensor $h$ related to the metric extensor $g$
(which determines $G$) and the relation between Hodge (star) extensors
associated to $G_{E}$ and $G$. \ In Section 4 we present our conclusions.

\section{ Geometric Algebras of Multivectors}

\subsection{Multivectors}

Let $V$ be $n$-dimensional vector space over the real field $\mathbb{R}$,
$T(V)$ the tensor algebra of $V$ and $\bigwedge V$ the exterior algebra of
$V$. The elements of $\bigwedge V$ are called multivectors.

We recall that $\bigwedge V$ is a $2^{n}$-dimensional \textit{associative}
algebra with unity. In addition, it is a $\mathbb{Z}$-graded algebra, i.e.,
\[
\bigwedge V=\bigoplus_{r=0}^{n}\bigwedge\nolimits^{r}V,
\]
and
\[
\bigwedge\nolimits^{r}V\wedge\bigwedge\nolimits^{s}V^{\ast}\subset
\bigwedge\nolimits^{r+s}V,
\]
$r,s\geq0$, where $\bigwedge^{r}V$ is the ${\binom{\;n\;}{\;r\;}}$-dimensional
subspace of \textit{homogenous} $r$-\textit{vectors\/} . $\ $We have moreover
the following identifications $\bigwedge^{0}V=\mathbb{R}$; $\bigwedge^{1}V=V$;
and of course, $\bigwedge^{r}V=\{0\}$ if $r>n$. \ If $A\in\bigwedge^{r}V$ for
some fixed $r$ ($r=0,\ldots,n$), then $A$ is said to be \textit{homogeneous}.
For any homogenous multivectors $A_{p}\in\bigwedge^{p}V$, $B_{q}\in
\bigwedge^{q}V$ we define here their \textit{exterior product}\footnote{If the
exterior algebra is defined as the quotient algebra $T(V)/J$, where $J$ is the
bilateral ideal generated by elements of the form $u\otimes v+v\otimes u$,
$u,v\in V$ then the exterior product $A_{p}\dot{\wedge}B_{q}$ is given by
$A_{p}\dot{\wedge}B_{q}=\mathcal{A}(A_{p}\otimes B_{q})$. Details on the
relation between $\wedge$ and $\overset{.}{\wedge}$ may be found in
\cite{10}.} by%
\begin{equation}
A_{p}\wedge B_{q}=\frac{(p+q)}{p!q!}\mathcal{A}(A_{p}\otimes B_{q}%
).\label{1.6biss}%
\end{equation}
where $\mathcal{A}:T_{k}(V\mathbf{)}\rightarrow\bigwedge^{k}V$ is the usual
\textit{antisymmetrization operator. }Of course, we have:
\begin{equation}
A_{p}\wedge B_{q}=(-1)^{pq}B_{p}\wedge A_{q}.
\end{equation}

For each $k=0,1,\ldots,n$ the linear mapping $\left\langle {}\right\rangle
_{k}:\bigwedge V\rightarrow\bigwedge^{k}V$ such that $X=X_{0}\oplus
X_{1}\oplus...\oplus X_{n}:=X_{0}+X_{1}+...+X_{n}$, $X_{k}\in\bigwedge^{k}V$,
$k=0,1,...,n$, then
\begin{equation}
\text{ }\left\langle X\right\rangle _{k}=X_{k}%
\end{equation}
is called the $k$\emph{-component projection operator} because $\left\langle
X\right\rangle _{k}$ is just the $k$-component of $X$.

For each $k=0,1,\ldots,n$ any multivector $X$ such that $\left\langle
X_{k}\right\rangle _{j}=0,$ for $k\neq j,$ is of course a \emph{homogeneous
multivector of degree }$k,$ or for short, a $k$\emph{-homogeneous
multivector}. It should be noticed that $0$ is an homogeneous multivector of
any degree $0,1,\ldots,n$.

Let $\{e_{j}\}$ be a basis of $V,$ and $\{\varepsilon^{j}\}$ be its dual basis
for $V^{\ast},$ i.e., $\varepsilon^{j}(e_{i})=\delta_{i}^{j}.$ Now, let us
take $t\in T^{k}V$ with $k\geq1.$ Such a \emph{contravariant }$k$%
-\emph{tensor} $t$ can be expanded onto the $k$-\emph{tensor basis}
$\{e_{j_{1}}\otimes\ldots\otimes e_{j_{k}}\}$ with $j_{1},\ldots
,j_{k}=1,\ldots,n$, by the well-known formula
\begin{equation}
t=t^{j_{1}\ldots j_{k}}e_{j_{1}}\otimes\ldots\otimes e_{j_{k}}, \label{EP.2e}%
\end{equation}
where $t^{j_{1}\ldots j_{k}}=t(\varepsilon^{j_{1}},\ldots,\varepsilon^{j_{k}%
})$ are the so-called $j_{1}\ldots j_{k}$-\emph{contravariant components} of
$t$ with respect to $\{e_{j_{1}}\otimes\ldots\otimes e_{j_{k}}\}.$

Using the definition of the antisymmetrization operator $\mathcal{A}$ it
follows (non-trivially) a remarkable identity which holds for the basis
$1$-forms $\varepsilon^{1},\ldots,\varepsilon^{n}$ belonging to $\{\varepsilon
^{j}\}.$ It is
\begin{equation}
\mathcal{A}t(\varepsilon^{j_{1}},\ldots,\varepsilon^{j_{k}})=\frac{1}%
{k!}\delta_{i_{1}\ldots i_{k}}^{j_{1}\ldots j_{k}}t(\varepsilon^{i_{1}}%
,\ldots,\varepsilon^{i_{k}}), \label{EP.2f}%
\end{equation}
where $\delta_{i_{1}\ldots i_{k}}^{j_{1}\ldots j_{k}}$ is the so-called
\emph{generalized Kronecker symbol of order} $k,$
\begin{equation}
\delta_{i_{1}\ldots i_{k}}^{j_{1}\ldots j_{k}}=\det\left[
\begin{array}
[c]{ccc}%
\delta_{i_{1}}^{j_{1}} & \ldots & \delta_{i_{1}}^{j_{k}}\\
\ldots & \ldots & \ldots\\
\delta_{i_{k}}^{j_{1}} & \ldots & \delta_{i_{k}}^{j_{k}}%
\end{array}
\right]  \text{ with }i_{1},\ldots,i_{k}\text{ and }j_{1},\ldots,j_{k}\text{
running from }1\text{ to }n. \label{EP.2g}%
\end{equation}

Let us take $X\in\bigwedge^{k}V$ with $k\geq2.$ By definition $X\in T^{k}V$
and is completely skew-symmetric, hence, it must be $X=\mathcal{A}X.$ Then, by
using Eq.(\ref{EP.2f}) we get a combinatorial identity which relates the
$i_{1}\ldots i_{k}$-components to the $j_{1}\ldots j_{k}$-components for $X.$
It is
\begin{equation}
X^{j_{1}\ldots j_{k}}=\frac{1}{k!}\delta_{i_{1}\ldots i_{k}}^{j_{1}\ldots
j_{k}}X^{i_{1}\ldots i_{k}}. \label{EP.2h}%
\end{equation}

From Eq.(\ref{1.6biss}) by using a well-known property of the
antisymmetrization operator, namely: $\mathcal{A}(\mathcal{A}t\otimes
u)=\mathcal{A}(t\otimes\mathcal{A}u)=\mathcal{A}(t\otimes u)$, we have the
following formula for expressing simple $k$-vectors in terms of the tensor
products of $k$ vectors. It is
\begin{equation}
v_{1}\wedge\ldots\wedge v_{k}=\epsilon^{i_{1}\ldots i_{k}}v_{i_{1}}%
\otimes\ldots\otimes v_{i_{k}}. \label{EP.3}%
\end{equation}
If $\omega^{1},\ldots,\omega^{k}\in V^{\ast},$ then
\begin{equation}
v_{1}\wedge\ldots\wedge v_{k}(\omega^{1},\ldots,\omega^{k})=\epsilon
^{i_{1}\ldots i_{k}}\omega^{1}(v_{i_{1}})\ldots\omega^{k}(v_{i_{k}}).
\label{EP.4}%
\end{equation}

Eq.(\ref{EP.3}) implies (non-trivially) a remarkable identity which holds for
the basis vectors $e_{1},\ldots,e_{n}$ belonging to any basis $\{e_{j}\}$ of
$V.$ It is
\begin{equation}
e_{i_{1}}\wedge\ldots\wedge e_{i_{k}}=\delta_{i_{1}\ldots i_{k}}^{j_{1}\ldots
j_{k}}e_{j_{1}}\otimes\ldots\otimes e_{j_{k}}. \label{EP.5a}%
\end{equation}

Once again let us take $X_{k}\in\bigwedge^{k}V$ with $k\geq2.$ Since $X_{k}\in
T^{k}V$ and is completely skew-symmetric, the use of Eq.(\ref{EP.2h}) and
Eq.(\ref{EP.5a}) in Eq.(\ref{EP.2e}) allows us to obtain the expansion
formula
\begin{equation}
X_{k}=\frac{1}{k!}X_{k}^{i_{1}\ldots i_{k}}e_{i_{1}}\wedge\ldots\wedge
e_{i_{k}}, \label{EP.5b}%
\end{equation}
where $X_{k}(\varepsilon^{j_{1}},\ldots,\varepsilon^{j_{k}})=X_{k}%
^{i_{1}\ldots i_{k}}$.

From this we see that a basis for the $2^{n}$-dimensional vector space of the
algebra $\bigwedge V$ $\ $is the set $\{e_{J}\}$ where $J$ are colective
indices, denoting specifically $\{1,e_{i},\frac{1}{2!}e_{i_{1}}\wedge
e_{i_{2}},...,\frac{1}{n!}e_{i_{1}}\wedge e_{i_{2}}...\wedge e_{i_{2}}\}$,
$i_{1},....i_{n}=0,1,2,...,n$. Then for a general multivector $X\in\bigwedge
V$ we have the expansion formula,%
\begin{equation}
X=s+v^{i}e_{i}+\frac{1}{2!}b^{ij}e_{i}\wedge e_{j}+\frac{1}{3!}t^{ijk}%
e_{i}\wedge e_{j}\wedge e_{k}+...+pe_{1}\wedge...\wedge e_{n}, \label{EP.5bis}%
\end{equation}
with $s,v^{i},b^{ij},t^{ijk},...,p\in\mathbb{R}$.

We recall moreover that the exterior product of $X,Y\in\bigwedge V,$ namely
$X\wedge Y\in\bigwedge V,$ is given by
\begin{equation}
X\wedge Y=\overset{n}{\underset{k=0}{\sum}}\overset{k}{\underset{j=0}{\sum}%
}\langle X\rangle_{j}\wedge\langle Y\rangle_{k-j}. \label{EP.6}%
\end{equation}

\subsection{Metric Structure}

Let us equip $V$ with a \emph{metric tensor}, $G:V\times V\rightarrow
\mathbb{R}$. As usual we write
\begin{equation}
G(v,w)\equiv v\cdot w, \label{MS.1c}%
\end{equation}
and call $v\cdot w$ the \emph{scalar product of the vectors }$v,w\in V.$

The pair $(V,G)$ is called a \emph{metric structure} for $V.$ Sometimes, $V$
is said to be a \emph{scalar product vector space}$.$

Let $\{e_{k}\}$ be any basis of $V,$ and $\{\varepsilon^{k}\}$ be its dual
basis for $V^{\ast}$. Let $G_{jk}=G(e_{j},e_{k})$, since $G$ is
non-degenerate, it follows that $\det\left[  G_{jk}\right]  \neq0.$ Then,
there exists the $jk$-entries for the inverse matrix of $\left[
G_{jk}\right]  ,$ namely $G^{jk},$ i.e., $G^{ks}G_{sj}=G_{js}G^{sk}=\delta
_{j}^{k}.$

We introduce the \emph{scalar product} \emph{of} $1$-\emph{forms}
$\omega,\sigma\in V^{\ast}$ by
\begin{equation}
\omega\cdot\sigma=G^{jk}\omega(e_{j})\sigma(e_{k}). \label{MS.2}%
\end{equation}
It should be noticed that the real number given by Eq.(\ref{MS.2}) does not
depend on the choice of $\{e_{k}\}$.

Now, we can define the so-called \textit{reciprocal bases }$\{e^{k}\}$ and
\textit{ }$\{\varepsilon_{k}\}$ of the bases $\{e_{k}\}$ and $\{\varepsilon
^{k}\}.$ Associated to $\{e_{k}\}$ we introduce the well-defined basis
$\{e^{k}\}$ by
\begin{equation}
e^{k}=G^{ks}e_{s},\text{ for each }k=1,\ldots,n. \label{MS.3a}%
\end{equation}
Such $e^{1},\ldots,e^{n}\in V$ are the unique basis vectors for $V$ which
satisfy
\begin{equation}
e^{k}\cdot e_{j}=\delta_{j}^{k}. \label{MS.3b}%
\end{equation}
Associated to $\{\varepsilon^{k}\},$ we can also introduce a well-defined
basis $\{\varepsilon_{k}\}$ by
\begin{equation}
\varepsilon_{k}=G_{ks}\varepsilon^{s},\text{ for each }k=1,\ldots,n.
\label{MS.4a}%
\end{equation}
Such $\varepsilon_{1},\ldots,\varepsilon_{n}\in V^{\ast}$ are the unique basis
$1$-forms for $V^{\ast}$ which satisfy
\begin{equation}
\varepsilon_{j}\cdot\varepsilon^{k}=\delta_{j}^{k}. \label{MS.4b}%
\end{equation}

The bases $\{e^{k}\}$ and $\{\varepsilon_{k}\}$ are respectively called the
reciprocal bases of $\{e_{k}\}$ and $\{\varepsilon^{k}\}$ (relative to the
metric tensor $G$).

Note that $\{\varepsilon_{k}\}$ is the dual basis of $\{e^{k}\},$ i.e.,
\begin{equation}
\varepsilon_{k}(e^{l})=\delta_{k}^{l}, \label{MS.5}%
\end{equation}
an immediate consequence of Eqs.(\ref{MS.4a}) and (\ref{MS.3a}).

From Eqs.(\ref{MS.3a}),( \ref{MS.3b}), (\ref{MS.4a}) and (\ref{MS.4b}), taking
into account Eq.(\ref{MS.2}), we easily get that
\begin{align}
\varepsilon_{j}\cdot\varepsilon_{k} &  =e_{j}\cdot e_{k},\label{MS.5a}\\
e^{j}\cdot e^{k} &  =G^{jk}=\varepsilon^{j}\cdot\varepsilon^{k}.\label{MS.5b}%
\end{align}

Using Eq.(\ref{MS.3b}) we get two expansion formulas for $v\in V$
\begin{equation}
v=v\cdot e^{k}e_{k}=v\cdot e_{k}e^{k}. \label{MS.6a}%
\end{equation}

Using Eq.(\ref{MS.4b}) we have that for all $\omega\in V^{\ast}$
\begin{equation}
\omega=\omega\cdot\varepsilon_{k}\varepsilon^{k}=\omega\cdot\varepsilon
^{k}\varepsilon_{k}. \label{MS.6b}%
\end{equation}

Let us take $X\in\bigwedge^{k}V$ with $k\geq2.$ By following analogous steps
to those which allowed us to get Eq.(\ref{EP.5b}) we can now obtain another
expansion formula for $k$-vectors, namely
\begin{equation}
X=\frac{1}{k!}X_{j_{1}\ldots j_{k}}e^{j_{1}}\wedge\ldots\wedge e^{j_{k}},
\label{MS.7}%
\end{equation}
where $X_{j_{1}\ldots j_{k}}=X(\varepsilon_{j_{1}},\ldots,\varepsilon_{j_{k}%
})$ are the so-called $j_{1}\ldots j_{k}$-\emph{covariant components of} $X$
(with respect to the $k$-tensor basis $\{e^{j_{1}}\otimes\ldots\otimes
e^{j_{k}}\}$ with $j_{1},\ldots,j_{k}=1,\ldots,n$).

Next, we will obtain a relation between the $i_{1}\ldots i_{k}$-covariant
components of $X$ and the $j_{1}\ldots j_{k}$-contravariant components of $X.$
A straightforward calculation yields
\begin{align*}
X(\varepsilon_{i_{1}}\ldots\varepsilon_{i_{k}})  &  =X(e_{i_{1}}\cdot
e_{s_{1}}\varepsilon^{s_{1}},\ldots,e_{i_{k}}\cdot e_{s_{k}}\varepsilon
^{s_{k}})\\
&  =X(\varepsilon^{s_{1}},\ldots,\varepsilon^{s_{k}})(e_{i_{1}}\cdot e_{s_{1}%
})\ldots(e_{i_{k}}\cdot e_{s_{k}})\\
&  =\frac{1}{k!}X(\varepsilon^{j_{1}},\ldots,\varepsilon^{j_{k}})\delta
_{j_{1}\ldots j_{k}}^{s_{1}\ldots s_{k}}(e_{i_{1}}\cdot e_{s_{1}}%
)\ldots(e_{i_{k}}\cdot e_{s_{k}}),
\end{align*}
hence,
\[
X_{i_{1}\ldots i_{k}}=\frac{1}{k!}X^{j_{1}\ldots j_{k}}\mathrm{Det}\left[
\begin{array}
[c]{ccc}%
e_{i_{1}}\cdot e_{j_{1}} & \ldots & e_{i_{1}}\cdot e_{j_{k}}\\
\ldots & \ldots & \ldots\\
e_{i_{k}}\cdot e_{j_{1}} & \ldots & e_{i_{k}}\cdot e_{j_{k}}%
\end{array}
\right]  .
\]

Finally we recall that we can take also as a basis the $2^{n}$-dimensional
vector space of the algebra $\bigwedge V$ the set $\{e^{J}\}$ where $J$ are
collective indices, denoting specifically $\{1,e^{i},\frac{1}{2!}e^{i_{1}%
}\wedge e^{_{i_{2}}},....,\frac{1}{n!}e^{i_{1}}\wedge e^{i_{2}}...e^{i_{n}}%
\}$, $i_{1},....i_{n}=0,1,2,...,n$. Then for a general multivector
$X\in\bigwedge V$ we have the expansion formula,%
\begin{equation}
X=s+v_{i}e^{i}+\frac{1}{2!}b_{ij}e^{i}\wedge e^{j}+\frac{1}{3!}t_{ijk}%
e^{i}\wedge e^{j}\wedge e^{k}+...+pe^{1}\wedge...\wedge e^{n},
\end{equation}
with $s,v_{i},b_{ij},t_{ijk},...,p\in\mathbb{R}$.

\subsection{Scalar Product for $\bigwedge^{p}V$}

Once a metric structure $(V,G)$ has been given we can equip $\bigwedge^{p}V$
with a scalar product of $p$-vectors. $\bigwedge V$ can then be endowed with a
scalar product of multivectors. This is done as follows.

The scalar product of $X_{p},Y_{p}\in\bigwedge^{p}V,$ namely $X_{p}\cdot
Y_{p}\in\mathbb{R},$ is defined by the axioms:

\textbf{Ax-i} For all $\alpha,\beta\in\mathbb{R}:$%
\begin{equation}
\alpha\cdot\beta=\alpha\beta\text{ (real product of }\alpha\text{ and }%
\beta\text{).} \label{SP.1a}%
\end{equation}

\textbf{Ax-ii} For all $X_{p},Y_{p}\in\bigwedge^{p}V,$ with $p\geq1:$%
\begin{align}
X_{p}\cdot Y_{p}  &  =\frac{1}{p!}X_{p}(\varepsilon^{i_{1}},\ldots
,\varepsilon^{i_{p}})Y_{p}(\varepsilon_{i_{1}},\ldots,\varepsilon_{i_{p}%
}),\nonumber\\
&  =\frac{1}{p!}X_{p}(\varepsilon_{i_{1}},\ldots,\varepsilon_{i_{p}}%
)Y_{p}(\varepsilon^{i_{1}},\ldots,\varepsilon^{i_{p}}), \label{SP.1b}%
\end{align}
where $\{\varepsilon_{i}\}$ is the reciprocal basis of $\{\varepsilon^{i}\},$
as defined by Eq.(\ref{MS.4a}).

It is a well-defined scalar product on $\bigwedge^{p}V,$ since it is
symmetric, satisfies the distributive laws, has the mixed associativity
property and is non-degenerate i.e., if $X_{p}\cdot Y_{p}=0$ for all $Y_{p},$
then $X_{p}=0$. For the special case of vectors Eq.(\ref{SP.1b}), of course,
reduces to
\begin{equation}
v\cdot w=\varepsilon^{i}(v)\varepsilon_{i}(w)=\varepsilon_{i}(v)\varepsilon
^{i}(w), \label{SP.1c}%
\end{equation}
i.e., $G=\varepsilon^{i}\otimes\varepsilon_{i}=\varepsilon_{i}\otimes
\varepsilon^{i}.$

The well-known formula for the scalar product of simple $k$-vectors can be
easily deduced from Eq.(\ref{SP.1b}). It is:
\begin{equation}
(v_{1}\wedge\ldots\wedge v_{k})\cdot(w_{1}\wedge\ldots\wedge w_{k}%
)=\mathrm{Det}\left[
\begin{array}
[c]{ccc}%
v_{1}\cdot w_{1} & \ldots & v_{1}\cdot w_{k}\\
\ldots & \ldots & \ldots\\
v_{k}\cdot w_{1} & \ldots & v_{k}\cdot w_{k}%
\end{array}
\right]  . \label{SP.2}%
\end{equation}

Now, we can generalize Eq.(\ref{MS.6a}) in order to get the expected expansion
formulas for $k$-vectors. For all $X\in\bigwedge^{k}V$ it holds two expansion
formulas
\begin{equation}
X=\frac{1}{k!}X\cdot(e^{j_{1}}\wedge\ldots e^{j_{k}})(e_{j_{1}}\wedge\ldots
e_{j_{k}})=\frac{1}{k!}X\cdot(e_{j_{1}}\wedge\ldots e_{j_{k}})(e^{j_{1}}%
\wedge\ldots e^{j_{k}}). \label{SP.2a}%
\end{equation}

and of course, $X^{j_{1}\ldots j_{k}}=X\cdot(e^{j_{1}}\wedge\ldots e^{j_{k}%
}).$ Analogously, we can prove that $X_{j_{1}\ldots j_{k}}=X\cdot(e_{j_{1}%
}\wedge\ldots e_{j_{k}})$.

\subsubsection{Scalar Product of Multivectors}

The scalar product of $X,Y\in\bigwedge V,$ namely $X\cdot Y\in\mathbb{R},$ is
defined by
\begin{equation}
X\cdot Y=\overset{n}{\underset{k=0}{\sum}}\langle X\rangle_{k}\cdot\langle
Y\rangle_{k}. \label{SP.3}%
\end{equation}

By using Eqs..(\ref{SP.1a}) and (\ref{SP.1b}) we can easily note that
Eq.(\ref{SP.3}) can still be written as
\begin{align}
X\cdot Y  &  =X_{0}Y_{0}+\overset{n}{\underset{k=1}{\sum}}\frac{1}{k!}%
X_{k}(\varepsilon^{i_{1}},\ldots,\varepsilon^{i_{k}})Y_{k}(\varepsilon_{i_{1}%
},\ldots,\varepsilon_{i_{k}})\nonumber\\
&  =X_{0}Y_{0}+\overset{n}{\underset{k=1}{\sum}}\frac{1}{k!}X_{k}%
(\varepsilon_{i_{1}},\ldots,\varepsilon_{i_{k}})Y_{k}(\varepsilon^{i_{1}%
},\ldots,\varepsilon^{i_{k}}). \label{SP.4}%
\end{align}

It is important to observe that the operation defined by Eq.(\ref{SP.3}) is
indeed a well-defined scalar product on $\bigwedge V$, since it is symmetric,
satisfies the distributive laws, has the mixed associative property and is not
degenerate, i.e., if $X\cdot Y=0$ for all $Y,$ then $X=0.$

\subsection{Involutions}

We recall that the main involution (or grade involution) denoted by $^{\wedge
}:%
%TCIMACRO{\dbigwedge }%
%BeginExpansion
{\displaystyle\bigwedge}
%EndExpansion
V\rightarrow%
%TCIMACRO{\dbigwedge }%
%BeginExpansion
{\displaystyle\bigwedge}
%EndExpansion
V$ satisfies: \ (i)\emph{ }if $\alpha\in\mathbb{R}$, $\hat{\alpha}=\alpha$;
(ii) if $a_{1}\wedge...\wedge a_{k}\in%
%TCIMACRO{\dbigwedge ^{k}}%
%BeginExpansion
{\displaystyle\bigwedge^{k}}
%EndExpansion
V$, $k\geq1$, $(a_{1}\wedge...\wedge a_{k})^{\symbol{94}}=(-1)^{k}a_{1}%
\wedge...\wedge a_{k}$; (iii) if $a,b\in\mathbb{R}$ and $\sigma,\tau\in%
%TCIMACRO{\dbigwedge ^{k}}%
%BeginExpansion
{\displaystyle\bigwedge^{k}}
%EndExpansion
V$ then $(a\sigma+b\tau)^{\ \wedge}=a\hat{\sigma}+b\hat{\tau}$; (iv) if $\tau=%
%TCIMACRO{\tsum }%
%BeginExpansion
{\textstyle\sum}
%EndExpansion
\tau_{k}$, $\tau_{k}\in%
%TCIMACRO{\dbigwedge ^{k}}%
%BeginExpansion
{\displaystyle\bigwedge^{k}}
%EndExpansion
V$ then%
\begin{equation}
\hat{\tau}=%
%TCIMACRO{\dsum \limits_{k}}%
%BeginExpansion
{\displaystyle\sum\limits_{k}}
%EndExpansion
\hat{\tau}_{k}. \label{T8}%
\end{equation}

We recall also that the reversion operator is the anti-automorphism
$\symbol{126}:%
%TCIMACRO{\dbigwedge }%
%BeginExpansion
{\displaystyle\bigwedge}
%EndExpansion
V\ni\mathbf{\tau\mapsto\tilde{\tau}\in}%
%TCIMACRO{\dbigwedge }%
%BeginExpansion
{\displaystyle\bigwedge}
%EndExpansion
V$ such that if $\tau=%
%TCIMACRO{\tsum }%
%BeginExpansion
{\textstyle\sum}
%EndExpansion
\langle\tau\rangle_{k}$, $\langle\tau\rangle_{k}\in%
%TCIMACRO{\dbigwedge ^{k}}%
%BeginExpansion
{\displaystyle\bigwedge^{k}}
%EndExpansion
V$ then:

(i)\emph{ }if $\alpha\in\mathbb{R}$, $\tilde{\alpha}=\alpha$; (ii) if
$a_{1}\wedge...\wedge a_{k}\in%
%TCIMACRO{\dbigwedge ^{k}}%
%BeginExpansion
{\displaystyle\bigwedge^{k}}
%EndExpansion
V$, $k\geq1$, $(a_{1}\wedge...\wedge a_{k})^{\symbol{126}}=a_{k}%
\wedge...\wedge a_{1}$; (iii) if $a,b\in\mathbb{R}$ and $\sigma,\tau\in%
%TCIMACRO{\dbigwedge ^{k}}%
%BeginExpansion
{\displaystyle\bigwedge^{k}}
%EndExpansion
V$ then $(a\sigma+b\tau)^{\symbol{126}}=a\tilde{\sigma}$ $+b\tilde{\tau}$;
(iv) if $\tau=%
%TCIMACRO{\tsum }%
%BeginExpansion
{\textstyle\sum}
%EndExpansion
\tau_{k}$, $\tau_{k}\in%
%TCIMACRO{\dbigwedge ^{k}}%
%BeginExpansion
{\displaystyle\bigwedge^{k}}
%EndExpansion
V$ then%
\begin{equation}
\tilde{\tau}=\sum\limits_{k=0}^{n}\tilde{\tau}_{k}, \label{T9}%
\end{equation}
where $\tilde{\tau}$\ is called the reverse of $\tau$.\medskip

Finally, we recall that the composition of the graded evolution with the
reversion operator, denoted by the symbol $-$ is called by some authors the
\emph{conjugation} and, $\bar{\tau}$ is said to be the \emph{conjugate} of
$\tau$. We have $\bar{\tau}=(\tilde{\tau})^{\symbol{94}}=(\hat{\tau})^{\sim}$.

\subsection{Contracted Products}

The left contracted product of $X_{p}\in\bigwedge^{p}V$ and $Y_{q}\in
\bigwedge^{q}V$ with $0\leq p\leq q\leq n,$ namely $X_{p}\lrcorner Y_{q}%
\in\bigwedge^{q-p}V,$ is defined by \ for all $X_{p}\in\bigwedge^{p}V$ and
$Y_{q}\in\bigwedge^{q}V$ with $p\leq q$ by:%
\begin{align}
X_{p}\lrcorner Y_{q}  &  =\frac{1}{(q-p)!}(\widetilde{X_{p}}\wedge e^{i_{1}%
}\wedge\ldots\wedge e^{i_{q-p}})\cdot Y_{q}e_{i_{1}}\wedge\ldots\wedge
e_{i_{q-p}}\nonumber\\
&  =\frac{1}{(q-p)!}(\widetilde{X_{p}}\wedge e_{i_{1}}\wedge\ldots\wedge
e_{i_{q-p}})\cdot Y_{q}e^{i_{1}}\wedge\ldots\wedge e^{i_{q-p}}. \label{CP.1b}%
\end{align}
It is clear that all $X_{p},Y_{p}\in\bigwedge^{p}V$ it holds%
\begin{equation}
X_{p}\lrcorner Y_{p}=\widetilde{X_{p}}\cdot Y_{p}=X_{p}\cdot\widetilde{Y_{p}}.
\label{CP.1a}%
\end{equation}
The right contracted product of $X_{p}\in\bigwedge^{p}V$ and $Y_{q}%
\in\bigwedge^{q}V$ with $n\geq p\geq q\geq0,$ namely $X_{p}\llcorner Y_{q}%
\in\bigwedge^{p-q}V,$ is defined for all $X_{p}\in\bigwedge^{p}V$ and
$Y_{q}\in\bigwedge^{q}V$ with $p>q$ by%
\begin{align}
X_{p}\llcorner Y_{q}  &  =\frac{1}{(p-q)!}X_{p}\cdot(e^{i_{1}}\wedge
\ldots\wedge e^{i_{p-q}}\wedge\widetilde{Y_{q}})e_{i_{1}}\wedge\ldots\wedge
e_{i_{p-q}}\nonumber\\
&  =\frac{1}{(p-q)!}X_{p}\cdot(e_{i_{1}}\wedge\ldots\wedge e_{i_{p-q}}%
\wedge\widetilde{Y_{q}})e^{i_{1}}\wedge\ldots\wedge e^{i_{p-q}}. \label{CP.2b}%
\end{align}
Of course, for all $X_{p},Y_{p}\in\bigwedge^{p}V$ we have%
\begin{equation}
X_{p}\llcorner Y_{p}=\widetilde{X_{p}}\cdot Y_{p}=X_{p}\cdot\widetilde{Y_{p}}.
\label{CP.2a}%
\end{equation}
It should be noticed that the $(q-p)$-vector defined by Eq.(\ref{CP.1b}) and
the $(p-q)$-vector defined by Eq.(\ref{CP.2b}) do not depend on the choice of
the reciprocal bases $\{e_{i}\}$ and $\{e^{i}\}$ used for calculating them.

Let us take $X_{p}\in\bigwedge^{p}V$ and $Y_{q}\in\bigwedge^{q}V$ with $p\leq
q.$ For all $Z_{q-p}\in\bigwedge^{q-p}V$ the following identity holds
\begin{equation}
(X_{p}\lrcorner Y_{q})\cdot Z_{q-p}=Y_{q}\cdot(\widetilde{X_{p}}\wedge
Z_{q-p}). \label{CP.3}%
\end{equation}
For $p<q$ Eq.(\ref{CP.3}) follows directly from Eq.(\ref{CP.1b}) and
Eq.(\ref{SP.2a}). But, for $p=q$ it trivially follows by taking into account
Eq.(\ref{CP.1a}), etc.

Let us take $X_{p}\in\bigwedge^{p}V$ and $Y_{q}\in\bigwedge^{q}V$ with $p\geq
q.$ For all $Z_{p-q}\in\bigwedge^{p-q}V$ the following identity holds
\begin{equation}
(X_{p}\llcorner Y_{q})\cdot Z_{p-q}=X_{p}\cdot(Z_{p-q}\wedge\widetilde{Y_{q}%
}). \label{CP.4}%
\end{equation}
For $p>q$ Eq.(\ref{CP.4}) follows directly from Eq.(\ref{CP.2b}) and
Eq.(\ref{SP.2a}). For $p=q$ it follows from Eq.(\ref{CP.2a}).

We recall moreover that for any $V_{p},W_{p}\in\bigwedge^{p}V$ and
$X_{q},Y_{q}\in\bigwedge^{q}V$ with $p\leq q$ we have $(V_{p}+W_{p})\lrcorner
X_{q}=V_{p}\lrcorner X_{q}+W_{p}\lrcorner X_{q}$, and $V_{p}\lrcorner
(X_{q}+Y_{q})=V_{p}\lrcorner X_{q}+V_{p}\lrcorner Y_{q}$ and also fo rany
$V_{p},W_{p}\in\bigwedge^{p}V$ and $X_{q},Y_{q}\in\bigwedge^{q}V$ with $p\geq
q$ we have $(V_{p}+W_{p})\llcorner X_{q}=V_{p}\llcorner X_{q}+W_{p}\llcorner
X_{q}$, and $V_{p}\llcorner(X_{q}+Y_{q})=V_{p}\llcorner X_{q}+V_{p}\llcorner
Y_{q}$. More important, we have for any $X_{p}\in\bigwedge^{p}V$ and $Y_{q}%
\in\bigwedge^{q}V$ with $p\leq q$%
\begin{equation}
X_{p}\lrcorner Y_{q}=(-1)^{p(q-p)}Y_{q}\llcorner X_{p}. \label{CP.12}%
\end{equation}
which follows by using Eq.(\ref{CP.3}), Eq.(\ref{CP.4}) and that $X_{p}\wedge
Y_{q}=(-1)^{pq}Y_{q}\wedge X_{p}.$ . Indeed, we have that $(X_{p}\lrcorner
Y_{q})\cdot Z_{q-p}=Y_{q}\cdot(\widetilde{X_{p}}\wedge Z_{q-p})=(-1)^{p(q-p)}%
Y_{q}\cdot(Z_{q-p}\wedge\widetilde{X_{p}})=(-1)^{p(q-p)}(Y_{q}\llcorner
X_{p})\cdot Z_{q-p}$, hence, by non-degeneracy of the scalar product, the
required result follows.

\subsubsection{Contracted Product of Nonhomogeneous Multivectors}

The left and right contracted products of $X,Y\in\bigwedge V,$ namely
$X\lrcorner Y\in\bigwedge V$ and $X\llcorner Y\in\bigwedge V,$ are defined by
\begin{align}
X\lrcorner Y  &  =\overset{n}{\underset{k=0}{\sum}}\overset{n-k}%
{\underset{j=0}{\sum}}\langle\langle X\rangle_{j}\lrcorner\langle
Y\rangle_{k+j}\rangle_{k}.\label{CP.5a}\\
X\llcorner Y  &  =\overset{n}{\underset{k=0}{\sum}}\overset{n-k}%
{\underset{j=0}{\sum}}\langle\langle X\rangle_{k+j}\llcorner\langle
Y\rangle_{j}\rangle_{k}. \label{CP.5b}%
\end{align}

We finalize this section presenting two noticeable formulas\footnote{For a
proof see \cite{10}} involving the contracted products and the scalar product,
and two other remarkable formulas relating the contracted products to the
exterior product and scalar product. They appear frequently in calculations

For any $X,Y,Z\in\bigwedge V$
\begin{align}
(X\lrcorner Y)\cdot Z  &  =Y\cdot(\widetilde{X}\wedge Z),\label{CP.8a}\\
(X\llcorner Y)\cdot Z  &  =X\cdot(Z\wedge\widetilde{Y}). \label{CP.8b}%
\end{align}

For any $X,Y,Z\in\bigwedge V$
\begin{align}
X\lrcorner(Y\lrcorner Z)  &  =(X\wedge Y)\cdot Z,\label{CP.9a}\\
(X\llcorner Y)\llcorner Z  &  =X\cdot(Y\wedge Z). \label{CP.9b}%
\end{align}

\subsection{Clifford Product and $\mathcal{C}\ell(V,G)$}

The two interior algebras together with the exterior algebra allow us to
define a \emph{Clifford product} of multivectors which is also an internal law
on $\bigwedge V.$ The Clifford product of $X,Y\in\bigwedge V,$ denoted by
juxtaposition $XY\in\bigwedge V,$ is defined by the following axioms

\textbf{Ax-ci} For all $\alpha\in\mathbb{R}$ and $X\in\bigwedge V:$%
\begin{equation}
\alpha X\text{ \ is the scalar multiplication of }X\text{ by }\alpha\text{.}
\label{CL.1a}%
\end{equation}

\textbf{Ax-cii} For all $v\in V$ and $X\in\bigwedge V:$%
\begin{align}
vX  &  =v\lrcorner X+v\wedge X,\label{CL.1b}\\
Xv  &  =X\llcorner v+X\wedge v. \label{CL.1c}%
\end{align}

\textbf{Ax-ciii} For all $X,Y,Z\in\bigwedge V:$%
\begin{equation}
(XY)Z=X(YZ). \label{CL.1d}%
\end{equation}

The Clifford product is distributive and associative. $\bigwedge V$ endowed
with this Clifford product is an associative algebra which will be called
the\emph{\ geometric algebra of multivectors} associated to a metric structure
$(V,G).$ It will be denoted by $\mathcal{C}\ell(V,G)$. Using the above axioms
we can derive a general formula for the Clifford\ product of two arbitrary
multivectors $A=%
%TCIMACRO{\dsum \nolimits_{r}}%
%BeginExpansion
{\displaystyle\sum\nolimits_{r}}
%EndExpansion
\oplus A_{r},B=%
%TCIMACRO{\dsum \nolimits_{s}}%
%BeginExpansion
{\displaystyle\sum\nolimits_{s}}
%EndExpansion
\oplus B_{s}\in\mathcal{C}\ell(V,G)$. We have%
\begin{align*}
AB &  =%
%TCIMACRO{\dsum \nolimits_{r,s}}%
%BeginExpansion
{\displaystyle\sum\nolimits_{r,s}}
%EndExpansion
\oplus A_{r}B_{s},\\
A_{r}B_{s} &  =\langle A_{r}B_{s}\rangle_{\left\vert r-s\right\vert }+\langle
A_{r}B_{s}\rangle_{\left\vert r-s+2\right\vert }+...+\langle A_{r}B_{s}%
\rangle_{r+s}%
\end{align*}

To continue we introduce one more convention. We denote by $X_{p}\ast Y_{q}$
either $(\wedge),$ or $(\cdot)$, or$(\lrcorner)$, or $(\llcorner)$ or
$($\emph{Clifford product}$).$

\subsection{ Euclidean and Metric Geometric Algebras}

\subsubsection{$\mathcal{C}\ell(V,G_{E})$}

Let us equip $V$ with an arbitrary (but fixed once for all) Euclidean metric
$G_{E}$. $V$ endowed with an Euclidean metric $G_{E},$ i.e., $(V,G_{E}),$ is
called an Euclidean metric structure for $V$. Sometimes, $(V,G_{E})$ is said
to be an Euclidean space.

Associated to $(V,G_{E})$ an Euclidean scalar product of vectors $v,w\in V$ is
given by
\begin{equation}
v\underset{G_{E}}{\cdot}w=G_{E}(v,w). \label{GA.2}%
\end{equation}
We introduce also an Euclidean scalar product of $p$-vectors $X_{p},Y_{p}%
\in\bigwedge^{p}V$ and Euclidean scalar product of multivectors $X,Y\in
\bigwedge V,$ namely $X_{p}\underset{G_{E}}{\cdot}Y_{p}\in\mathbb{R}$ and
$X\underset{G_{E}}{\cdot}Y\in\mathbb{R},$ using respectively the
Eqs.(\ref{SP.1a}) and (\ref{SP.1b}), and Eq.(\ref{SP.3}). The Clifford algebra
associated to the pair $(V,G_{E})$ will be denoted $\mathcal{C}\ell(V,G_{E})$
and called \textit{Euclidean geometric algebra}. It will play a role in our
theory a role analogous to the one of matrix algebra in standard presentations
of linear algebra.

\subsubsection{$\mathcal{C}\ell(V,G)$}

Let us take any metric tensor $G$ on the vector space $V.$ Associated to the
metric structure $(V,G)$ a scalar product of vectors $v,w\in V$ is represented
by
\begin{equation}
v\underset{G}{\cdot}w=G(v,w). \label{GA.3}%
\end{equation}
Of course, the corresponding scalar product of $p$-vectors $X_{p},Y_{p}%
\in\bigwedge^{p}V$ and scalar product of multivectors $X,Y\in\bigwedge V,$
namely $X_{p}\underset{G}{\cdot}Y_{p}\in\mathbb{R}$ and $X\underset{G}{\cdot
}Y\in\mathbb{R},$ are defined respectively by Eqs.(\ref{SP.1a}) and
(\ref{SP.1b}) and Eq.(\ref{SP.3}). The Clifford algebra associated to the pair
$(V,G)$ will be denoted $\mathcal{C}\ell(V,G)$ and called \textit{metric
geometric algebra.}

We will find a relationship between $(V,G)$ and $(V,G_{E})$, thereby showing
how an arbitrary $G$-scalar product on $\bigwedge^{p}V$ and $\bigwedge V$ is
related to a $G_{E}$-scalar products on $\bigwedge^{p}V$ and $\bigwedge V$.
This starting point which permits us to relate $\mathcal{C}\ell(V,G)$ with
$\mathcal{C}\ell(V,G_{E})$ is the concept of a metric operator $g$ (but a
convenient algorithm needs the concept of deformation extensors to be
introduced in Section 3) which we now introduce.

\subsubsection{Enter $g$}

To continue, choose once and for all a \textit{fiducial} Euclidean metric
structure $(V,G_{E})$. We now recall that for any metric tensor $G$ there
exists an unique linear operator $g:V\rightarrow V$, such that for all $v,w\in
V$
\begin{equation}
v\underset{G}{\cdot}w=g(v)\underset{G_{E}}{\cdot}w.\label{GA.4}%
\end{equation}
Such $g$ is given by
\begin{equation}
g(v)=(v\underset{G}{\cdot}e_{k})\underset{G_{E}}{e^{k}}=(v\underset{G}{\cdot
}\underset{G_{E}}{e^{k}})e_{k},\label{GA.5}%
\end{equation}
where $\{e_{k}\}$ is any basis of $V,$ and $\{\underset{G_{E}}{e^{k}}\}$ is
its reciprocal basis with respect to $(V,G_{E}),$ i.e., $e_{k}\underset{G_{E}%
}{\cdot}\underset{G_{E}}{e^{l}}=\delta_{k}^{l}.$ Note that the vector $g(v)$
does not depend on the basis $\{e_{k}\}$ chosen for calculating it.

We now show that $g(v)$ given by Eq.(\ref{GA.5}) satisfies Eq.(\ref{GA.4}).
Using Eq.(\ref{MS.6a}) we have,
\[
g(v)\underset{G_{E}}{\cdot}w=(v\underset{G}{\cdot}e_{k})(\underset{G_{E}%
}{e^{k}}\underset{G_{E}}{\cdot}w)=(v\underset{G}{\cdot}(\underset{G_{E}}%
{e^{k}}\underset{G_{E}}{\cdot}w)e_{k})=v\underset{G}{\cdot}w.
\]

Now, suppose that there is some $g^{\prime}$ which satisfies Eq.(\ref{GA.4}),
i.e., $v\underset{G}{\cdot}w=g^{\prime}(v)\underset{G_{E}}{\cdot}w.$ Then,
using once again Eq.(\ref{MS.6a}) we have
\[
g^{\prime}(v)=(g^{\prime}(v)\underset{G_{E}}{\cdot}e_{k})\underset{G_{E}%
}{e^{k}}=(v\underset{G}{\cdot}e_{k})\underset{G_{E}}{e^{k}}=g(v),
\]
i.e., $g^{\prime}=g$. So the existence and the uniqueness of such a linear
operator $g$ are proved.

Since $G$ is a symmetric covariant $2$-tensor over $V$, i.e., $G(v,w)=G(w,v)$
$\forall v,w\in V$, it follows from Eq.(\ref{GA.4}) that $g$ is an
\emph{adjoint symmetric} linear operator with respect to $(V,G_{E})$, i.e.,
\begin{equation}
g(v)\underset{G_{E}}{\cdot}w=v\underset{G_{E}}{\cdot}g(w), \label{GA.6}%
\end{equation}
The property expressed by Eq.(\ref{GA.6}) may be coded as $g=g^{\dagger
(G_{E})}$.

Since $G$ is a non-degenerate covariant $2$-tensor over $V$ (i.e., if
$G(v,w)=0$ $\forall w\in V,$ then $v=0$) it follows that $g$ is a non-singular
(invertible) linear operator. Its inverse linear operator is given by the
noticeable formula
\begin{equation}
g^{-1}(v)=G^{jk}(v\underset{G_{E}}{\cdot}e_{j})e_{k}, \label{GA.7}%
\end{equation}
where $G^{jk}$ are the $jk$-entries of the inverse matrix of $\left[
G_{jk}\right]  $ with $G_{jk}\equiv G(e_{j},e_{k}).$ Note that the vector
$g^{-1}(v) $ does not depend on the basis $\{e_{k}\}$ chosen for its calculation.

We must prove that indeed $g^{-1}\circ g=g\circ g^{-1}=i_{V},$ where $i_{V}$
is the identity function for $V.$

By using Eq.(\ref{GA.7}), Eq.(\ref{GA.4}), Eq.(\ref{MS.3a}) for $(V,G)$ and
Eq.(\ref{MS.6a}) for $(V,G),$ we immediately have that $g^{-1}\circ g(v)=v$,
i.e., $g^{-1}\circ g=i_{V}$. Also, $g\circ g^{-1}(v)=v$, i.e., $g\circ
g^{-1}=i_{V}$.

It should be remarked that such $g$ only depends on the \textit{choice} of the
fiducial Euclidean structure $(V,G_{E}).$ However, $g$ codifies all the
geometric information contained in $G.$ Such $g$ will be called the
\emph{metric operator} for $G.$

Now, we show how the scalar product $X_{p}\underset{G}{\cdot}Y_{p}$ is related
to the scalar product $X_{p}\underset{G_{E}}{\cdot}Y_{p}.$

For any simple $k$-vectors $v_{1}\wedge\ldots\wedge v_{k}\in\bigwedge^{k}V$
and $w_{1}\wedge\ldots\wedge w_{k}\in\bigwedge^{k}V$ it holds
\begin{equation}
(v_{1}\wedge\ldots v_{k})\underset{G}{\cdot}(w_{1}\wedge\ldots w_{k}%
)=(g(v_{1})\wedge\ldots\wedge g(v_{k}))\underset{G_{E}}{\cdot}(w_{1}%
\wedge\ldots w_{k}), \label{GA.8}%
\end{equation}

To show Eq.(\ref{GA.8}) we will use Eq.(\ref{SP.2}) for $(V,G)$ and
$(V,G_{E}).$ By using Eq.(\ref{GA.4}) and a fundamental property of any
outermorphism, a straightforward calculation yields
\[
(v_{1}\wedge\ldots v_{k})\underset{G}{\cdot}(w_{1}\wedge\ldots w_{k}%
)=(g(v_{1})\wedge\ldots\wedge g(v_{k}))\underset{G_{E}}{\cdot}(w_{1}%
\wedge\ldots w_{k})
\]

For any $k$-vectors $X_{k},Y_{k}\in\bigwedge^{k}V$ it holds
\begin{equation}
X_{k}\underset{G}{\cdot}Y_{k}=\underline{g}(X_{k})\underset{G_{E}}{\cdot}%
Y_{k}. \label{GA.9}%
\end{equation}
\medskip where if $X_{k}=\frac{1}{k!}X_{i_{1}...i_{k}}e_{i_{1}}\wedge
\ldots\wedge e_{i_{k}}$, $\underline{g}(X_{k})=\frac{1}{k!}X_{i_{1}...i_{k}%
}g(e_{i_{1}})\wedge\ldots\wedge g(e_{i_{k}})$. The operator $\underline{g}$ is
called the exterior power extension of $g$ and the general properties of a
general exterior power operator are discussed in Section 3.3..

Next, we show how the scalar product $X\underset{G}{\cdot}Y$ is related to the
scalar product $X\underset{G_{E}}{\cdot}Y.$

We also recall that for any multivectors $X,Y\in\bigwedge V$ it holds
\begin{equation}
X\underset{G}{\cdot}Y=\underline{g}(X)\underset{G_{E}}{\cdot}Y. \label{GA.10}%
\end{equation}
\medskip

The $G$-contracted products are related to the $G_{E}$-contracted products by
two noticeable formulas.

For any $X,Y\in\bigwedge V$
\begin{align}
X\underset{G}{\lrcorner}Y  &  =\underline{g}(X)\underset{G_{E}}{\lrcorner
}Y,\label{GA.11a}\\
X\underset{G}{\llcorner}Y  &  =X\underset{G_{E}}{\llcorner}\underline{g}(Y).
\label{GA.11b}%
\end{align}

The Clifford algebra associated to the pair $(V,G)$ will be denoted
$\mathcal{C}\ell(V,G)$ or $\mathcal{C}\ell(V,g)$ and called \textit{metric
geometric algebra}. We shall use also, in what follows, the notation
\ $X\underset{g}{\cdot}Y$ meaning  $X\underset{G}{\cdot}Y$.

\section{ Theory of Extensors}

In this section we recall some basic notions of the theory of extensors thus
completing the presentation of the algebraic notions necessary for the
remaining papers of the series.\footnote{We recall that extensors are a
particular case of a more important concept, namely of multivector functions
of multivector variables. These objects and still more general ones called
multivector functionals are used, e.g., in the Lagrangian formulation of the
theory of\ multivector and extensor fields. Their theory may be found in
\cite{10,rodoliv2006}.}

\subsection{General $k$-Extensors}

Let $\bigwedge_{1}^{\diamond}V,\ldots,$ $\bigwedge_{k}^{\diamond}V$ be $k$
subspaces\footnote{We can have, e.g., $\bigwedge_{1}^{\diamond}V=%
%TCIMACRO{\dbigwedge \nolimits^{3}}%
%BeginExpansion
{\displaystyle\bigwedge\nolimits^{3}}
%EndExpansion
V\oplus%
%TCIMACRO{\dbigwedge \nolimits^{5}}%
%BeginExpansion
{\displaystyle\bigwedge\nolimits^{5}}
%EndExpansion
V$, $\bigwedge_{2}^{\diamond}V=%
%TCIMACRO{\dbigwedge \nolimits^{1}}%
%BeginExpansion
{\displaystyle\bigwedge\nolimits^{1}}
%EndExpansion
V\oplus%
%TCIMACRO{\dbigwedge \nolimits^{3}}%
%BeginExpansion
{\displaystyle\bigwedge\nolimits^{3}}
%EndExpansion
V\oplus%
%TCIMACRO{\dbigwedge \nolimits^{n}}%
%BeginExpansion
{\displaystyle\bigwedge\nolimits^{n}}
%EndExpansion
$, etc.} of $\bigwedge V$ such that each of them is \emph{any} sum of
homogeneous subspaces of $\bigwedge V$, and $\bigwedge^{\diamond}V$ is either
any sum of homogeneous subspaces of $\bigwedge V$ or even the trivial subspace
consisting of the null vector $\{0\}$. A multilinear mapping from the
cartesian product $\bigwedge_{1}^{\diamond}V\times\cdots\times\bigwedge
_{k}^{\diamond}V$ to $\bigwedge^{\diamond}V$ will be called a general
$k$-extensor over $V,$ i.e., $t:\bigwedge_{1}^{\diamond}V\times\cdots
\times\bigwedge_{k}^{\diamond}V\rightarrow\bigwedge^{\diamond}V$ such that for
any $\alpha_{j},\alpha_{j}^{\prime}\in\mathbb{R}$ and $X_{j},X_{j}^{\prime}%
\in\bigwedge_{j}^{\diamond}V,$
\begin{equation}
t(\ldots,\alpha_{j}X_{j}+\alpha_{j}^{\prime}X_{j}^{\prime},\ldots)=\alpha
_{j}t(\ldots,X_{j},\ldots)+\alpha_{j}^{\prime}t(\ldots,X_{j}^{\prime},\ldots),
\label{GkE.1}%
\end{equation}
for each $j$ with $1\leq j\leq k.$

It should be noticed that the linear operators on $V,$ $\bigwedge^{p}V$ or
$\bigwedge V$ which appear in ordinary linear algebra are particular cases of
$1$-extensors over $V.$ Note also that a covariant $k$-tensor over $V$ is just
a $k$-extensor over $V.$ On this way, the concept of general $k$-extensor
generalizes and unifies both of the concepts of linear operator and of
covariant $k$-tensor. These mathematical objects are of the same nature!

The set of general $k$-extensors over $V,$ denoted by $k$-$ext(\bigwedge
_{1}^{\diamond}V,\ldots,\bigwedge_{k}^{\diamond}V;\bigwedge^{\diamond}V),$ has
a natural structure of real vector space. Its dimension is clearly given by
\begin{equation}
\dim k\text{-}ext(%
%TCIMACRO{\dbigwedge \nolimits_{1}^{\diamond}}%
%BeginExpansion
{\displaystyle\bigwedge\nolimits_{1}^{\diamond}}
%EndExpansion
V,\ldots,%
%TCIMACRO{\dbigwedge \nolimits_{k}^{\diamond}}%
%BeginExpansion
{\displaystyle\bigwedge\nolimits_{k}^{\diamond}}
%EndExpansion
V;%
%TCIMACRO{\dbigwedge \nolimits^{\diamond}}%
%BeginExpansion
{\displaystyle\bigwedge\nolimits^{\diamond}}
%EndExpansion
V)=\dim%
%TCIMACRO{\dbigwedge \nolimits_{1}^{\diamond}}%
%BeginExpansion
{\displaystyle\bigwedge\nolimits_{1}^{\diamond}}
%EndExpansion
V\cdots\dim%
%TCIMACRO{\dbigwedge \nolimits_{k}^{\diamond}}%
%BeginExpansion
{\displaystyle\bigwedge\nolimits_{k}^{\diamond}}
%EndExpansion
V\dim%
%TCIMACRO{\dbigwedge \nolimits^{\diamond}}%
%BeginExpansion
{\displaystyle\bigwedge\nolimits^{\diamond}}
%EndExpansion
V. \label{GE.2}%
\end{equation}
We shall need to consider only some particular cases of these general
$k$-extensors over $V.$ So, special names and notations will be given for them.

We will equip $V$ with an arbitrary (but fixed once and for all) Euclidean
metric $G_{E}$, and denote the scalar product of multivectors $X,Y\in\bigwedge
V$ with respect to the Euclidean metric structure $(V,G_{E}),$ $X\cdot Y$
instead of the more detailed notation $X\underset{G_{E}}{\cdot}Y$.

\subsubsection{$(p,q)$-Extensors}

Let $\{e_{j}\}$ be any basis for $V,$ and $\{e^{j}\}$ be its Euclidean
reciprocal basis for $V,$ i.e., $e_{j}\cdot e^{k}=\delta_{j}^{k}$. Let $p$ and
$q$ be two integer numbers with $0\leq p,q\leq n$. A linear mapping which
sends $p$-vectors to $q$-vectors will be called a $(p,q)$-\emph{extensor} over
$V.$ The space of these objects, namely $1$-$ext(\bigwedge^{p}V;\bigwedge
^{q}V),$ will be denoted by $ext_{p}^{q}(V)$ for short. By using
Eq.(\ref{GE.2}) we get
\begin{equation}
\dim ext_{p}^{q}(V)=\binom{n}{p}\binom{n}{q}. \label{GE.3'}%
\end{equation}

For instance, we see that the $(1,1)$-extensors over $V$ are just the
well-known linear operators on $V$.

To continue we recall that an Euclidean scalar product $\cdot$ has been
intoduced in order to be possible to define the Euclidean Clifford Algebra
$\mathcal{C\ell(}V,G_{E})$ which, as we hope it is clear at this point, is the
basic instrument in our calculations. This \textit{suggests} the introduction
of the following basis for the space of extensors.

Let $\varepsilon^{j_{1}\ldots j_{p};k_{1}\ldots k_{q}}\in ext_{p}^{q}(V)$ be
$\binom{n}{p}\binom{n}{q}$ extensors such that
\begin{equation}
\varepsilon^{j_{1}\ldots j_{p};k_{1}\ldots k_{q}}(X)=(e^{j_{1}}\wedge
\ldots\wedge e^{j_{p}})\cdot Xe^{k_{1}}\wedge\ldots\wedge e^{k_{q}%
}.\label{GE.3aa}%
\end{equation}
We show now that they define a a $(p,q)$-extensor basis for $ext_{p}^{q}(V).$

Indeed, the extensors given by Eq.(\ref{GE.3aa}) are linearly independent, and
for each $t\in ext_{p}^{q}(V)$ there exist $\binom{n}{p}\binom{n}{q}$ real
numbers, say $t_{j_{1}\ldots j_{p};k_{1}\ldots k_{q}},$ given by
\begin{equation}
t_{j_{1}\ldots j_{p};k_{1}\ldots k_{q}}=t(e_{j_{1}}\wedge\ldots\wedge
e_{j_{p}})\cdot(e_{k_{1}}\wedge\ldots\wedge e_{k_{q}}) \label{GE.3bb}%
\end{equation}
such that
\begin{equation}
t=\frac{1}{p!}\frac{1}{q!}t_{j_{1}\ldots j_{p};k_{1}\ldots k_{q}}%
\varepsilon^{j_{1}\ldots j_{p};k_{1}\ldots k_{q}}. \label{GE.3cc}%
\end{equation}
Such $t_{j_{1}\ldots j_{p};k_{1}\ldots k_{q}}$ will be called the $j_{1}\ldots
j_{p};k_{1}\ldots k_{q}$-th \emph{covariant components} of $t$ with respect to
the $(p,q)$-\emph{extensor basis} $\{\varepsilon^{j_{1}\ldots j_{p}%
;k_{1}\ldots k_{q}}\}.$

Of course, there are still other kinds of $(p,q)$-extensor bases for
$ext_{p}^{q}(V)$ besides the one given by Eq.(\ref{GE.3aa}) which can be
constructed from the vector bases $\{e_{j}\}$ and $\{e^{j}\}.$ The total
number of these different kinds of $(p,q)$-extensor bases for $ext_{p}^{q}(V)$
are $2^{p+q}$.

Now, if we take the basis $(p,q)$-extensors $\varepsilon_{j_{1}\ldots
j_{p};k_{1}\ldots k_{q}}$ and the real numbers $t^{j_{1}\ldots j_{p}%
;k_{1}\ldots k_{q}}$ defined by
\begin{align}
\varepsilon_{j_{1}\ldots j_{p};k_{1}\ldots k_{q}}(X)  &  =(e_{j_{1}}%
\wedge\ldots\wedge e_{j_{p}})\cdot Xe_{k_{1}}\wedge\ldots\wedge e_{k_{q}%
},\label{GE.3d}\\
t^{j_{1}\ldots j_{p};k_{1}\ldots k_{q}}  &  =t(e^{j_{1}}\wedge\ldots\wedge
e^{j_{p}})\cdot(e^{k_{1}}\wedge\ldots\wedge e^{k_{q}}), \label{GE.3e}%
\end{align}
we get an expansion formula for $t\in ext_{p}^{q}(V)$ analogous to that given
by Eq.(\ref{GE.3c}), i.e.,
\begin{equation}
t=\frac{1}{p!}\frac{1}{q!}t^{j_{1}\ldots j_{p};k_{1}\ldots k_{q}}%
\varepsilon_{j_{1}\ldots j_{p};k_{1}\ldots k_{q}}. \label{GE.3f}%
\end{equation}
Such $t^{j_{1}\ldots j_{p};k_{1}\ldots k_{q}}$ are called the $j_{1}\ldots
j_{p};k_{1}\ldots k_{q}$-th \emph{contravariant components} of $t$ with
respect to the $(p,q)$-\emph{extensor basis} $\{\varepsilon_{j_{1}\ldots
j_{p};k_{1}\ldots k_{q}}\}.$

\subsubsection{ Extensors}

A linear mapping which sends multivectors to multivectors will be simply
called an \emph{extensor} over $V.$ They are the linear operators on
$\bigwedge V.$ For the space of extensors over $V,$ namely $1$-$ext(\bigwedge
V;\bigwedge V),$ we will use the short notation $ext(V).$ By using
Eq.(\ref{GE.2}) we get
\begin{equation}
\dim ext(V)=2^{n}2^{n}. \label{GE.4}%
\end{equation}

For instance, we will see that the so-called Hodge star operator is just a
well-defined extensor over $V$ which can be thought as an exterior direct sum
of $(p,n-p)$-extensor over $V$. \footnote{The extended (or exterior power) of
$t\in ext_{1}^{1}(V)$ as defined in Section 3.3. is just an extensor over $V,$
i.e., $\underline{t}\in ext(V).$}

There are $2^{n}2^{n}$ extensors over $V,$ namely $\varepsilon^{J;K}$, given
by\footnote{Recall once again that $J$ and $K$ are colective indices,
$e_{J}=1,e_{j_{1}},e_{j_{1}}\wedge e_{j_{2}},\ldots$($e^{J}=1,e^{j_{1}%
},e^{j_{1}}\wedge e^{j_{2}},\ldots$).}
\begin{equation}
\varepsilon^{J;K}(X)=(e^{J}\cdot X)e^{K} \label{GE.4a}%
\end{equation}
which can be used to introduce an \emph{\ extensor} \emph{basis }for $ext(V).$

In fact they are linearly independent, and for each $t\in ext(V)$ there exist
$2^{n}2^{n}$ real numbers, say $t_{J;K},$ given by
\begin{equation}
t_{J;K}=t(e_{J})\cdot e_{K} \label{GE.4b}%
\end{equation}
such that
\begin{equation}
t=\underset{J}{\sum}\underset{K}{\sum}\frac{1}{\nu(J)!}\frac{1}{\nu
(K)!}t_{J;K}\varepsilon^{J;K}, \label{GE.4c}%
\end{equation}
where we define $\nu(J)=0,1,2,\ldots$ for $J=\emptyset$, $j_{1},j_{1}%
j_{2},\ldots,$ where all index $j_{1},j_{2},\ldots$ runs from $1$ to $n$. Such
$t_{J;K}$ will be called the $J;K$-th\emph{\ covariant components} of $t$ with
respect to the \emph{extensor bases} $\{\varepsilon^{J;K}\}.$

We notice that exactly $(2^{n+1}-1)^{2}$ extensor bases for $ext(V)$ can be
constructed from the basis vectors $\{e_{j}\}$ and $\{e^{j}\}.$ For instance,
whenever the basis extensors $\varepsilon_{J;K}$ and the real numbers
$t^{J;K}$ defined by
\begin{align}
\varepsilon_{J;K}(X)  &  =(e_{J}\cdot X)e_{K},\label{GE.4d}\\
t^{J;K}  &  =t(e^{J})\cdot e^{K} \label{GE.4e}%
\end{align}
are used, an expansion formula for $t\in ext(V)$ analogous to that given by
Eq.(\ref{GE.4c}) can be obtained, i.e.,
\begin{equation}
t=\underset{J}{\sum}\underset{K}{\sum}\frac{1}{\nu(J)!}\frac{1}{\nu
(K)!}t^{J;K}\varepsilon_{J;K}. \label{GE.4f}%
\end{equation}
Such $t^{J;K}$ are called the $J;K$-th \emph{contravariant components} of $t$
with respect to the \emph{extensor bases} $\{\varepsilon_{J;K}\}$.

\subsubsection{Elementary $k$-Extensors}

A multilinear mapping which takes $k$-uple of vectors into $q$-vectors will be
called an \emph{elementary }$k$-\emph{extensor} over $V$ of degree $q.$ The
space of these objects, namely $k$-$ext(V,\ldots,V;\bigwedge^{q}V),$ will be
denoted by $k$-$ext^{q}(V).$ It is easy to verify (using Eq.(\ref{GE.2}))
that
\begin{equation}
\dim k\text{-}ext^{q}(V)=n^{k}\binom{n}{q}. \label{GE.5}%
\end{equation}
It should be noticed that an\emph{\ }elementary $k$-extensor over $V$ of
degree $0$ is just a \emph{covariant }$k$-\emph{tensor} over $V,$ i.e.,
$k$-$ext^{0}(V)\equiv T_{k}(V).$ It is easily realized that $1$-$ext^{q}%
(V)\equiv ext_{1}^{q}(V).$

The elementary $k$-extensors of degrees $0,1,2,\ldots$  are sometimes said to
be \emph{scalar}, \emph{vector,} \emph{bivector}, $\ldots$ \emph{elementary
}$k$-\emph{extensors.}

The $n^{k}\binom{n}{q}$ elementary $k$-extensors of degree $q$ belonging to
$k$-$ext^{q}(V),$ namely $\varepsilon^{j_{1},\ldots,j_{k};k_{1}\ldots k_{q}},$
given by
\begin{equation}
\varepsilon^{j_{1},\ldots,j_{k};k_{1}\ldots k_{q}}(v_{1},\ldots,v_{k}%
)=(v_{1}\cdot e^{j_{1}})\ldots(v_{k}\cdot e^{j_{k}})e^{k_{1}}\wedge
\ldots\wedge e^{k_{q}} \label{GE.5a}%
\end{equation}
define elementary \emph{basis vectors}, ( i.e., $k$-extensor of degree $q$)
for $k$-$ext^{q}(V)$.

In fact they are linearly independent, and for all $t\in k$-$ext^{q}(V)$ there
are $n^{k}\binom{n}{q}$ real numbers, say $t_{j_{1},\ldots,j_{k};k_{1}\ldots
k_{q}},$ given by
\begin{equation}
t_{j_{1},\ldots,j_{k};k_{1}\ldots k_{q}}=t(e_{j_{1}},\ldots,e_{j_{k}}%
)\cdot(e_{k_{1}}\wedge\ldots\wedge e_{k_{q}}) \label{GE.5b}%
\end{equation}
such that
\begin{equation}
t=\frac{1}{q!}t_{j_{1},\ldots,j_{k};k_{1}\ldots k_{q}}\varepsilon
^{j_{1},\ldots,j_{k};k_{1}\ldots k_{q}}. \label{GE.5c}%
\end{equation}
Such $t_{j_{1},\ldots,j_{k};k_{1}\ldots k_{q}}$ will be called the
$j_{1},\ldots,j_{k};k_{1}\ldots k_{q}$-th \emph{covariant components }of $t$
with respect to the\emph{\ basis} $\{\varepsilon^{j_{1},\ldots,j_{k}%
;k_{1}\ldots k_{q}}\}.$

We notice that exactly $2^{k+q}$ elementary $k$-extensors of degree $q$ bases
for $k$-$ext^{q}(V)$ can be constructed from the vector bases $\{e_{j}\}$ and
$\{e^{j}\}.$ For instance, we may define $\varepsilon_{j_{1},\ldots
,j_{k};k_{1}\ldots k_{q}}$ (the basis elementary $k$-extensor of degree $q$)
and the real numbers $t^{j_{1},\ldots,j_{k};k_{1}\ldots k_{q}}$ by
\begin{align}
\varepsilon_{j_{1},\ldots,j_{k};k_{1}\ldots k_{q}}(v_{1},\ldots,v_{k})  &
=(v_{1}\cdot e_{j_{1}})\ldots(v_{k}\cdot e_{j_{k}})e_{k_{1}}\wedge\ldots\wedge
e_{k_{q}},\label{GE.5d}\\
t^{j_{1},\ldots,j_{k};k_{1}\ldots k_{q}}  &  =t(e^{j_{1}},\ldots,e^{j_{k}%
})\cdot(e^{k_{1}}\wedge\ldots\wedge e^{k_{q}}). \label{GE.5e}%
\end{align}
Then, we also have other expansion formulas for $t\in k$-$ext^{q}(V)$ besides
that given by Eq.(\ref{GE.5c}), e.g.,
\begin{equation}
t=\frac{1}{q!}t^{j_{1},\ldots,j_{k};k_{1}\ldots k_{q}}\varepsilon
_{j_{1},\ldots,j_{k};k_{1}\ldots k_{q}}. \label{GE.5f}%
\end{equation}
Such $t^{j_{1},\ldots,j_{k};k_{1}\ldots k_{q}}$ are called the $j_{1}%
,\ldots,j_{k};k_{1}\ldots k_{q}$-th \emph{contravariant components} of $t$
with respect to the \emph{basis} $\{\varepsilon_{j_{1},\ldots,j_{k}%
;k_{1}\ldots k_{q}}\}$.

\subsection{Projectors}

Let $\bigwedge^{\diamond}V$ be either any sum of homogeneous
subspaces\footnote{Note that for such a subspace $\bigwedge^{\diamond}V$ there
are $\nu$ integers $p_{1,}\ldots,p_{\nu}$ ($0\leq p_{1}<\cdots<p_{\nu}\leq n$)
such that $\bigwedge_{1}^{\diamond}V=\bigwedge^{p_{1}}V\oplus\cdots
\oplus\bigwedge^{p_{\nu}}V$.} of $\bigwedge V$ or the trivial subspace
$\{0\}.$ Associated to $\bigwedge^{\diamond}V,$ a noticeable extensor from
$\bigwedge V$ to $\bigwedge^{\diamond}V,$ namely $\left\langle \left.
{}\right.  \right\rangle _{\bigwedge^{\diamond}V},$ can defined by
\begin{equation}
\left\langle X\right\rangle _{\Lambda^{\diamond}V}=\left\{
\begin{array}
[c]{cc}%
\left\langle X\right\rangle _{p_{1}}+\cdots+\left\langle X\right\rangle
_{p_{v}}, & \text{if }\bigwedge^{\diamond}V=\bigwedge^{p_{1}}V\oplus
\cdots\oplus\bigwedge^{p_{\nu}}V\\
0, & \text{if }\bigwedge^{\diamond}V=\{0\}
\end{array}
\right.  . \label{P.1}%
\end{equation}
Such $\left\langle \left.  {}\right.  \right\rangle _{\bigwedge^{\diamond}%
V}\in1$-$ext(\bigwedge V;\bigwedge^{\diamond}V)$ will be called the
$\bigwedge^{\diamond}V$-\emph{projector extensor}.

We notice that if $\bigwedge^{\diamond}V$ is any homogeneous subspace of
$\bigwedge V,$ i.e., $\bigwedge^{\diamond}V=\bigwedge^{p}V,$ then the
projector extensor is reduced to the so-called $p$-\emph{part operator}, i.e.,
$\left\langle \left.  {}\right.  \right\rangle _{\bigwedge^{\diamond}%
V}=\left\langle \left.  {}\right.  \right\rangle _{p}.$

We now summarize the fundamental properties for the $\bigwedge^{\diamond}%
V$-projector extensors.

Let $\bigwedge_{1}^{\diamond}V$ and $\bigwedge_{2}^{\diamond}V$ be two
subspaces of $\bigwedge V.$ If each of them is either any sum of homogeneous
subspaces of $\bigwedge V$ or the trivial subspace $\{0\},$ then
\begin{align}
\left\langle \left\langle X\right\rangle _{\bigwedge_{1}^{\diamond}%
V}\right\rangle _{\bigwedge_{2}^{\diamond}V}  &  =\left\langle X\right\rangle
_{\bigwedge_{1}^{\diamond}V\cap\bigwedge_{2}^{\diamond}V}\label{P.2a}\\
\left\langle X\right\rangle _{\bigwedge_{1}^{\diamond}V}+\left\langle
X\right\rangle _{\bigwedge_{2}^{\diamond}V}  &  =\left\langle X\right\rangle
_{\bigwedge_{1}^{\diamond}V\cup\bigwedge_{2}^{\diamond}V}. \label{P.2b}%
\end{align}

Let $\bigwedge^{\diamond}V$ be either any sum of homogeneous subspaces of
$\bigwedge V$ or the trivial subspace $\{0\}.$ Then, it holds
\begin{equation}
\left\langle X\right\rangle _{\bigwedge^{\diamond}V}\cdot Y=X\cdot\left\langle
Y\right\rangle _{\bigwedge^{\diamond}V}.. \label{P.2c}%
\end{equation}

We see that the concept of $\bigwedge^{\diamond}V$-projector extensor is just
a natural generalization of the concept of $p$-part operator.

\subsection{ Exterior Power Extension Operator}

Let $\{e_{j}\}$ be any basis for $V,$ and $\{\varepsilon^{j}\}$ be its dual
basis for $V^{\ast}.$ As we know, $\{\varepsilon^{j}\}$ is the unique $1$-form
basis associated to the vector basis $\{e_{j}\}$ such that $\varepsilon
^{j}(e_{i})=\delta_{i}^{j}.$ The linear mapping $ext_{1}^{1}(V)\ni
t\mapsto\underline{t}\in ext(V)$ such that for any $X\in\bigwedge V$ and
$X=X_{0}+\overset{n}{\underset{k=1}{%
%TCIMACRO{\dsum }%
%BeginExpansion
{\displaystyle\sum}
%EndExpansion
}}X_{k},$ then
\begin{equation}
\underline{t}(X)=X_{0}+\overset{n}{\underset{k=1}{%
%TCIMACRO{\dsum }%
%BeginExpansion
{\displaystyle\sum}
%EndExpansion
}}\frac{1}{k!}X_{k}(\varepsilon^{j_{1}},\ldots,\varepsilon^{j_{k}})t(e_{j_{1}%
})\wedge\ldots\wedge t(e_{j_{k}}) \label{EO.1}%
\end{equation}
will be called the exterior power \emph{extension operator}, or extension
operator for short. We call $\underline{t}$ the \emph{extended} of $t.$ It is
the well-known outermorphism of $t$ in ordinary linear algebra.

The extension operator is well-defined since it does not depend on the choice
of $\{e_{j}\}.$

We summarize now the basic properties satisfied by the extension operator.

\textbf{e1} The extension operator is grade-preserving, i.e.,
\begin{equation}
\text{if }X\in\bigwedge^{p}V,\text{ then }\underline{t}(X)\in\bigwedge^{p}V.
\label{EO.1a}%
\end{equation}
It is an obvious result which follows from Eq.(\ref{EO.1}).

\textbf{e2} For any $\alpha\in\mathbb{R},$ $v\in V$ and $v_{1}\wedge
\ldots\wedge v_{k}\in\bigwedge^{k}V$,
\begin{align}
\underline{t}(\alpha)  &  =\alpha,\label{EO.2a}\\
\underline{t}(v)  &  =t(v),\label{EO.2b}\\
\underline{t}(v_{1}\wedge\ldots\wedge v_{k})  &  =t(v_{1})\wedge\ldots\wedge
t(v_{k}). \label{EO.2c}%
\end{align}

\textbf{e3} For any $X,Y\in\bigwedge V$,
\begin{equation}
\underline{t}(X\wedge Y)=\underline{t}(X)\wedge\underline{t}(Y). \label{EO.3}%
\end{equation}
Eq.(\ref{EO.3}) an immediate result which follows from Eq.(\ref{EO.2c}).

We emphasize that the three fundamental properties as given by Eq.(\ref{EO.2a}%
), Eq.(\ref{EO.2b}) and Eq.(\ref{EO.3}) together are completely equivalent to
the \emph{extension procedure} as defined by Eq.(\ref{EO.1}).

We present next some important properties for the extension operator.

\textbf{e4} Let us take $s,t\in ext_{1}^{1}(V).$ Then, the following result
holds
\begin{equation}
\underline{s\circ t}=\underline{s}\circ\underline{t}. \label{EO.4}%
\end{equation}

It is enough to present the proofs for scalars and simple $k$-vectors.

For $\alpha\in\mathbb{R},$ by using Eq.(\ref{EO.2a}) we get
\[
\underline{s\circ t}(\alpha)=\alpha=\underline{s}(\alpha)=\underline
{s}(\underline{t}(\alpha))=\underline{s}\circ\underline{t}(\alpha).
\]

For a simple $k$-vector $v_{1}\wedge\ldots\wedge v_{k}\in\bigwedge^{k}V,$ by
using Eq.(\ref{EO.2c}) we get
\begin{align*}
\underline{s\circ t}(v_{1}\wedge\ldots\wedge v_{k})  &  =s\circ t(v_{1}%
)\wedge\ldots\wedge s\circ t(v_{k})=s(t(v_{1}))\wedge\ldots\wedge
s(t(v_{k}))\\
&  =\underline{s}(t(v_{1})\wedge\ldots\wedge t(v_{k}))=\underline
{s}(\underline{t}(v_{1}\wedge\ldots\wedge v_{k})),\\
&  =\underline{s}\circ\underline{t}(v_{1}\wedge\ldots\wedge v_{k}).
\end{align*}

Next, we can easily generalize to multivectors due to the linearity of
extensors. It yields
\[
\underline{s\circ t}(X)=\underline{s}\circ\underline{t}(X).
\]

\textbf{e5} Let us take $t\in ext_{1}^{1}(V)$ with inverse $t^{-1}\in
ext_{1}^{1}(V),$ i.e., $t^{-1}\circ t=t\circ t^{-1}=i_{V}.$ Then,
$\underline{(t^{-1})}\in ext(V)$ is the inverse of $\underline{t}\in ext(V),$
i.e.,
\begin{equation}
(\underline{t})^{-1}=\underline{(t^{-1})}. \label{EO.5}%
\end{equation}
Indeed, by using Eq.(\ref{EO.4}) and the obvious property $\underline{i}%
_{V}=i_{\bigwedge V},$ we have that
\[
t^{-1}\circ t=t\circ t^{-1}=i_{V}\Rightarrow\underline{(t^{-1})}%
\circ\underline{t}=\underline{t}\circ\underline{(t^{-1})}=i_{\bigwedge V},
\]
which means that the inverse of the extended of $t$ equals the extended of the
inverse of $t.$

In accordance with the above corollary we use in what follows a more simple
notation $\underline{t}^{-1}$ to denote both $(\underline{t})^{-1}$ and
$\underline{(t^{-1})}.$

Let $\{e_{j}\}$ be any basis for $V,$ and $\{e^{j}\}$ its Euclidean reciprocal
basis for $V,$ i.e., $e_{j}\cdot e^{k}=\delta_{j}^{k}.$ There are two
interesting and useful formulas for calculating the extended of $t\in
ext_{1}^{1}(V),$ i.e.,%

\begin{align}
\underline{t}(X)  &  =1\cdot X+\overset{n}{\underset{k=1}{\sum}}\frac{1}%
{k!}(e^{j_{1}}\wedge\ldots\wedge e^{j_{k}})\cdot Xt(e_{j_{1}})\wedge
\ldots\wedge t(e_{j_{k}})\label{EO.6a}\\
&  =1\cdot X+\overset{n}{\underset{k=1}{\sum}}\frac{1}{k!}(e_{j_{1}}%
\wedge\ldots\wedge e_{j_{k}})\cdot Xt(e^{j_{1}})\wedge\ldots\wedge t(e^{j_{k}%
}). \label{EO.6b}%
\end{align}

\subsection{Standard Adjoint Operator}

Let $\bigwedge_{1}^{\diamond}V$ and $\bigwedge_{2}^{\diamond}V$ be two
subspaces of $\bigwedge V$ such that each of them is any sum of homogeneous
subspaces of $\bigwedge V.$ Let $\{e_{j}\}$ and $\{e^{j}\}$ be two Euclidean
reciprocal bases to each other for $V,$ i.e., $e_{j}\cdot e^{k}=\delta_{j}%
^{k}.$

We call \emph{standard adjoint operator }of $t$ the linear mapping
$1$-$ext(\bigwedge_{1}^{\diamond}V;\bigwedge_{2}^{\diamond}V)\ni t\rightarrow
t^{\dagger}\in1$-$ext(\bigwedge_{2}^{\diamond}V;\bigwedge_{1}^{\diamond}V)$
such that for any $Y\in\bigwedge_{2}^{\diamond}V:$%
\begin{align}
t^{\dagger}(Y)  &  =t(\left\langle 1\right\rangle _{\bigwedge_{1}^{\diamond}%
V})\cdot Y+\underset{k=1}{\overset{n}{\sum}}\frac{1}{k!}t(\left\langle
e^{j_{1}}\wedge\ldots e^{j_{k}}\right\rangle _{\bigwedge_{1}^{\diamond}%
V})\cdot Ye_{j_{1}}\wedge\ldots e_{j_{k}}\label{SAO.1a}\\
&  =t(\left\langle 1\right\rangle _{\bigwedge_{1}^{\diamond}V})\cdot
Y+\underset{k=1}{\overset{n}{\sum}}\frac{1}{k!}t(\left\langle e_{j_{1}}%
\wedge\ldots e_{j_{k}}\right\rangle _{\bigwedge_{1}^{\diamond}V})\cdot
Ye^{j_{1}}\wedge\ldots e^{j_{k}}. \label{SAO.1b}%
\end{align}
Using a more compact notation by employing the \emph{collective index} $J$ we
can write
\begin{align}
t^{\dagger}(Y)  &  =\underset{J}{\sum}\frac{1}{\nu(J)!}t(\left\langle
e^{J}\right\rangle _{\bigwedge_{1}^{\diamond}V})\cdot Ye_{J}\label{SAO.2a}\\
&  =\underset{J}{\sum}\frac{1}{\nu(J)!}t(\left\langle e_{J}\right\rangle
_{\bigwedge_{1}^{\diamond}V})\cdot Ye^{J}, \label{SAO.2b}%
\end{align}
We call $t^{\dagger}$ the \emph{standard adjoint} of $t.$ It should be noticed
the use of the $\bigwedge_{1}^{\diamond}V$-projector extensor.

The standard adjoint operator is well-defined since the sums appearing in each
one of the above places do not depend on the choice of $\{e_{j}\}.$

Let us take $X\in\bigwedge_{1}^{\diamond}V$ and $Y\in\bigwedge_{2}^{\diamond
}V.$ A straightforward calculation yields
\begin{align*}
X\cdot t^{\dagger}(Y)  &  =\underset{J}{\sum}\frac{1}{\nu(J)!}t(\left\langle
e^{J}\right\rangle _{\bigwedge_{1}^{\diamond}V})\cdot Y(X\cdot e_{J})\\
&  =t(\underset{J}{\sum}\frac{1}{\nu(J)!}\left\langle (X\cdot e_{J}%
)e^{J}\right\rangle _{\bigwedge_{1}^{\diamond}V})\cdot Y\\
&  =t(\left\langle X\right\rangle _{\bigwedge_{1}^{\diamond}V})\cdot Y,
\end{align*}
i.e.,
\begin{equation}
X\cdot t^{\dagger}(Y)=t(X)\cdot Y. \label{SAO.3}%
\end{equation}
It is a generalization of the well-known property which holds for linear operators.

Let us take $t\in1$-$ext(\bigwedge_{1}^{\diamond}V;\bigwedge_{2}^{\diamond}V)$
and $u\in1$-$ext(\bigwedge_{2}^{\diamond}V;\bigwedge_{3}^{\diamond}V).$ We can
note that $u\circ t\in1$-$ext(\bigwedge_{1}^{\diamond}V;\bigwedge
_{3}^{\diamond}V)$ and $t^{\dagger}\circ u^{\dagger}\in1$-$ext(\bigwedge
_{3}^{\diamond}V;\bigwedge_{1}^{\diamond}V).$ Then, let us take $X\in
\bigwedge_{1}^{\diamond}V$ and $Z\in\bigwedge_{3}^{\diamond}V,$ by using
Eq.(\ref{SAO.3}) we have that
\[
X\cdot(u\circ t)^{\dagger}(Z)=(u\circ t)(X)\cdot Z=t(X)\cdot u^{\dagger
}(Z)=X\cdot(t^{\dagger}\circ u^{\dagger})(Z).
\]
Hence, we get
\begin{equation}
(u\circ t)^{\dagger}=t^{\dagger}\circ u^{\dagger}. \label{SAO.4}%
\end{equation}

Let us take $t\in1$-$ext(\bigwedge^{\diamond}V;\bigwedge^{\diamond}V)$ with
inverse $t^{-1}\in1$-$ext(\bigwedge^{\diamond}V;\bigwedge^{\diamond}V),$ i.e.,
$t^{-1}\circ t=t\circ t^{-1}=i_{\bigwedge^{\diamond}V},$ where $i_{\bigwedge
^{\diamond}V}\in1$-$ext(\bigwedge^{\diamond}V;\bigwedge^{\diamond}V)$ is the
so-called identity function for $\bigwedge^{\diamond}V.$ By using
Eq.(\ref{SAO.4}) and the obvious property $i_{\bigwedge^{\diamond}%
V}=i_{\bigwedge^{\diamond}V}^{\dagger},$ we have that
\[
t^{-1}\circ t=t\circ t^{-1}=i_{\bigwedge^{\diamond}V}\Rightarrow t^{\dagger
}\circ(t^{-1})^{\dagger}=(t^{-1})^{\dagger}\circ t^{\dagger}=i_{\bigwedge
^{\diamond}V},
\]
hence,
\begin{equation}
(t^{\dagger})^{-1}=(t^{-1})^{\dagger}, \label{SAO.5}%
\end{equation}
i.e., the inverse of the adjoint of $t$ equals the adjoint of the inverse of
$t$. In accordance with the above corollary it is possible to use a more
simple symbol, say $t^{\ast}$, to denote both of $(t^{\dagger})^{-1}$ and
$(t^{-1})^{\dagger}.$

Let us take $t\in ext_{1}^{1}(V).$ We note that $\underline{t}\in ext(V)$ and
$\underline{(t^{\dagger})}\in ext(V).$ A straightforward calculation by using
Eqs.(\ref{EO.6a}) and (\ref{EO.6b}) yields
\begin{align*}
\underline{(t^{\dagger})}(Y)  &  =1\cdot Y+\overset{n}{\underset{k=1}{\sum}%
}\frac{1}{k!}(e^{j_{1}}\wedge\ldots e^{j_{k}})\cdot Yt^{\dagger}(e_{j_{1}%
})\wedge\ldots t^{\dagger}(e_{j_{k}})\\
&  =1\cdot Y+\\
&  \overset{n}{\underset{k=1}{\sum}}\frac{1}{k!}(e^{j_{1}}\wedge\ldots
e^{j_{k}})\cdot Yt^{\dagger}(e_{j_{1}})\cdot e_{p_{1}}e^{p_{1}}\wedge\ldots
t^{\dagger}(e_{j_{k}})\cdot e_{p_{k}}e^{p_{k}}\\
&  =1\cdot Y+\overset{n}{\underset{k=1}{\sum}}\frac{1}{k!}(e_{j_{1}}\cdot
t(e_{p_{1}})e^{j_{1}}\wedge\ldots e_{j_{k}}\cdot t(e_{p_{k}})e^{j_{k}})\cdot
Ye^{p_{1}}\wedge\ldots e^{p_{k}}\\
&  =\underline{t}(1)\cdot Y+\overset{n}{\underset{k=1}{\sum}}\frac{1}%
{k!}\underline{t}(e_{p_{1}}\wedge\ldots e_{p_{k}})\cdot Ye^{p_{1}}\wedge\ldots
e^{p_{k}}\\
&  =(\underline{t})^{\dagger}(Y).
\end{align*}
Hence, we get
\begin{equation}
\underline{(t^{\dagger})}=(\underline{t})^{\dagger}. \label{SAO.6}%
\end{equation}
This means that the extension operator commutes with the adjoint operator. In
accordance with the above property we may use a more simple notation
$\underline{t}^{\dagger}$ to denote without ambiguities both $\underline
{(t^{\dagger})}$ and $(\underline{t})^{\dagger}.$

\subsection{Standard Generalization Operator}

Let $\{e_{k}\}$ be any basis for $V,$ and $\{e^{k}\}$ be its Euclidean
reciprocal basis for $V,$ as we know, $e_{k}\cdot e^{l}=\delta_{k}^{l}.$

The linear mapping $ext_{1}^{1}(V)\ni t\mapsto\underset{\thicksim}{t}\in
ext(V)$ such that for any $X\in\bigwedge V$
\begin{equation}
\underset{\thicksim}{t}(X)=t(e^{k})\wedge(e_{k}\lrcorner X)=t(e_{k}%
)\wedge(e^{k}\lrcorner X) \label{SGO.1}%
\end{equation}
will be called the \emph{generalization operator.} We call $\underset
{\thicksim}{t}$ the \emph{generalized of }$t.$

The generalization operator is well-defined since it does not depend on the
choice of $\{e_{k}\}.$

We present now some important properties which are satisfied by the
generalization operator.\medskip

\textbf{g1 }The generalization operator is grade-preserving, i.e.,
\begin{equation}
\text{if }X\in\bigwedge^{k}V,\text{ then }\underset{\thicksim}{t}%
(X)\in\bigwedge^{k}V. \label{SGO.1a}%
\end{equation}
\medskip

\textbf{g2 }The grade involution $\widehat{\left.  {}\right.  }\in ext(V),$
reversion $\widetilde{\left.  {}\right.  }\in ext(V),$ and conjugation
$\overline{\left.  {}\right.  }\in ext(V)$ commute with the generalization
operator, i.e.,
\begin{align}
\underset{\thicksim}{t}(\widehat{X})  &  =\widehat{\underset{\thicksim}{t}%
(X)},\label{SGO.2a}\\
\underset{\thicksim}{t}(\widetilde{X})  &  =\widetilde{\underset{\thicksim}%
{t}(X)},\label{SGO.2b}\\
\underset{\thicksim}{t}(\overline{X})  &  =\overline{\underset{\thicksim}%
{t}(X)}. \label{SGO.2c}%
\end{align}
They are immediate consequences of the grade-preserving property.\medskip

\textbf{g3 }For any $\alpha\in\mathbb{R},$ $v\in V$ and $X,Y\in\bigwedge V $
it holds
\begin{align}
\underset{\thicksim}{t}(\alpha)  &  =0,\label{SGO.3a}\\
\underset{\thicksim}{t}(v)  &  =t(v),\label{SGO.3b}\\
\underset{\thicksim}{t}(X\wedge Y)  &  =\underset{\thicksim}{t}(X)\wedge
Y+X\wedge\underset{\thicksim}{t}(Y). \label{SGO.3c}%
\end{align}

We can show that the basic properties given by Eq.(\ref{SGO.3a}),
Eq.(\ref{SGO.3b}) and Eq.(\ref{SGO.3c}) together are completely equivalent to
the \emph{generalization procedure} as defined by Eq.(\ref{SGO.1}).\medskip

\textbf{g4} The generalization operator commutes with the adjoint operator,
i.e.,
\begin{equation}
(\underset{\thicksim}{t})^{\dagger}=\underset{\thicksim}{(t^{\dagger})},
\label{SGO.4}%
\end{equation}
or put it on another way, the adjoint of the generalized of $t$ is just the
generalized of the adjoint of $t$.

The proof of this result is a straightforward calculation which uses
Eq.(\ref{SAO.3}) and the multivector identities: $X\cdot(a\wedge
Y)=(a\lrcorner X)\wedge Y$ and $X\cdot(a\lrcorner Y)=(a\wedge X)\cdot Y,$ with
$a\in V$ and $X,Y\in\bigwedge V$. Indeed,
\begin{align*}
(\underset{\thicksim}{t})^{\dagger}(X)\cdot Y  &  =X\cdot\underset{\thicksim
}{t}(Y)\\
&  =(e_{j}\wedge(t(e^{j})\lrcorner X))\cdot Y=(e_{j}\wedge(t(e^{j})\cdot
e^{k}e_{k}\lrcorner X))\cdot Y\\
&  =(e^{j}\cdot t^{\dagger}(e^{k})e_{j}\wedge(e_{k}\lrcorner X))\cdot
Y=(t^{\dagger}(e^{k})\wedge(e_{k}\lrcorner X))\cdot Y\\
&  =\underset{\thicksim}{(t^{\dagger})}(X)\cdot Y.
\end{align*}
Hence, by the non-degeneracy property of the Euclidean scalar product, the
required result follows.

In agreement with the above property we use in what follows a more simple
symbol, $\underset{\thicksim}{t^{\dagger}}$ to denote both $(\underset
{\thicksim}{t})^{\dagger}$ or $\underset{\thicksim}{(t^{\dagger})}.$

\textbf{g5} The symmetric (skew-symmetric) part of the generalized of $t$ is
just the generalized of the symmetric (skew-symmetric) part of $t,$ i.e.,
\begin{equation}
(\underset{\thicksim}{t})_{\pm}=(\underset{\thicksim}{t_{\pm}}). \label{SGO.5}%
\end{equation}
This property follows immediately from Eq.(\ref{SGO.4}).

We see also that it is possible to use a more simple notation, $\underset
{\thicksim}{t}_{\pm}$to denote $(\underset{\thicksim}{t})_{\pm}$ or
$(\underset{\thicksim}{t_{\pm}}).$

\textbf{g6} The skew-symmetric part of the generalized of $t$ can be
factorized by the noticeable formula\footnote{Recall that $X\times
Y\equiv\frac{1}{2}(XY-YX).$}
\begin{equation}
\underset{\thicksim}{t}_{-}(X)=\frac{1}{2}biv[t]\times X, \label{SGO.6}%
\end{equation}
where $biv[t]\equiv t(e^{k})\wedge e_{k}$ is a \emph{characteristic invariant}
of $t$, called the \emph{bivector of} $t.\medskip$

We prove this result recalling Eq.(\ref{SGO.5}), the well-known identity
$t_{-}(a)=\frac{1}{2}biv[t]\times a$ and the multivector identity $B\times
X=(B\times e^{k})\wedge(e_{k}\lrcorner X),$ with $B\in\bigwedge^{2}V$ and
$X\in\bigwedge V$. We have that
\[
\underset{\thicksim}{t}_{-}(X)=t_{-}(e^{k})\wedge(e_{k}\lrcorner X)=(\frac
{1}{2}biv[t]\times e^{k})\wedge(e_{k}\lrcorner X)=\frac{1}{2}biv[t]\times X.
\]

\textbf{g7} A noticeable formula holds for the skew-symmetric part of the
generalized of $t.$ For all $X,Y\in\bigwedge V$
\begin{equation}
\underset{\thicksim}{t}_{-}(X\ast Y)=\underset{\thicksim}{t}_{-}(X)\ast
Y+X\ast\underset{\thicksim}{t}_{-}(Y), \label{SGO.7}%
\end{equation}
where $\ast$ is any product either $(\wedge),$ $(\cdot),$ $(\lrcorner
,\llcorner)$ or $($\emph{Clifford product}$).$

In order to prove this property we must use Eq.(\ref{SGO.6}) and the
multivector identity $B\times(X\ast Y)=(B\times X)\ast Y+X\ast(B\times Y),$
with $B\in\bigwedge^{2}V$ and $X,Y\in\bigwedge V$. By taking into account
Eq.(\ref{SGO.3a}) we can see that the following property for the Euclidean
scalar product of multivectors holds
\begin{equation}
\underset{\thicksim}{t}_{-}(X)\cdot Y+X\cdot\underset{\thicksim}{t}_{-}(Y)=0.
\label{SGO.7a}%
\end{equation}
It is consistent with the well-known property: \emph{the adjoint of a
skew-symmetric extensor equals minus the extensor!}

\subsection{Determinant}

We define the \emph{determinant}\footnote{The concept of determinant of a
$(1,1)$-extensor is related, but distinct from the well known determinant of a
square matrix. For details the reader is invited to consult \cite{rodoliv2006}%
.}\emph{\ of }$t\in ext_{1}^{1}(V)$ as the mapping $\det:ext_{1}^{1}(V)\ni
t\mapsto\det[t]\in%
%TCIMACRO{\dbigwedge \nolimits^{0}}%
%BeginExpansion
{\displaystyle\bigwedge\nolimits^{0}}
%EndExpansion
V\equiv\mathbb{R},$ such that for all non-zero pseudoscalar $I\in%
%TCIMACRO{\dbigwedge \nolimits^{n}}%
%BeginExpansion
{\displaystyle\bigwedge\nolimits^{n}}
%EndExpansion
V$
\begin{equation}
\underline{t}(I)=\det[t]I.\label{D.1}%
\end{equation}

As can easily the value of $\det[t]$ does not depend\emph{ }on the choice of
the pseudo-scalar $I$. We recall that we choose the symbol $\det[t]$ for the
determinant of a $(1,1)$-extensor in order to not confuse this concept with
the concept of the determinant of a square matrix (see below).

We present now some of the most important properties satisfied by the
$\det[t]$.

\textbf{d1} Let $t$ and $u$ be two $(1,1)$--extensors. It holds
\begin{equation}
\det[u\circ t]=\det[u]\det[t]. \label{D.2}%
\end{equation}

Indeed, take a non-zero pseudoscalar $I\in\bigwedge^{n}V.$ By using
Eq.(\ref{EO.4}) and Eq.(\ref{D.1}) we can write that
\begin{align*}
\det[u\circ t]I  &  =\underline{u\circ t}(I)=\underline{u}\circ\underline
{t}(I)=\underline{u}(\underline{t}(I))\\
&  =\underline{u}(\det[t]I)=\det[t]\underline{u}(I),\\
&  =\det[t]\det[u]I.\blacksquare
\end{align*}

\textbf{d2} Let us take $t\in ext_{1}^{1}(V)$ with inverse $t^{-1}\in
ext_{1}^{1}(V).$ It holds
\begin{equation}
\det[t^{-1}]=(\det[t])^{-1}. \label{D.3}%
\end{equation}

Indeed, by using Eq.(\ref{D.2}) and the obvious property $\det[i_{V}]=1,$ we
have that
\[
t^{-1}\circ t=t\circ t^{-1}=i_{V}\Rightarrow\det[t^{-1}]\det[t]=\det
[t]\det[t^{-1}]=1,
\]
which means that the \emph{determinant of the inverse equals the inverse of
the determinant.}

Due to the above corollary it is often convenient to use the short notation
$\left.  \det\right.  ^{-1}[t]$ for both $\det[t^{-1}]$ and $(\det[t])^{-1}. $

\textbf{d3} Let us take $t\in ext_{1}^{1}(V).$ It holds
\begin{equation}
\det[t^{\dagger}]=\det[t]. \label{D.4}%
\end{equation}
Indeed, take a non-zero $I\in\bigwedge^{n}V.$ Then, by using Eq.(\ref{D.1})
and Eq.(\ref{SAO.3}) we have that
\[
\det[t^{\dagger}]I\cdot I=\underline{t}^{\dagger}(I)\cdot I=I\cdot
\underline{t}(I)=I\cdot\det[t]I=\det[t]I\cdot I,
\]
whence, the expected result follows.

Let $\{e_{j}\}$ be any basis for $V,$ and $\{e^{j}\}$ be its Euclidean
reciprocal basis for $V,$ i.e., $e_{j}\cdot e^{k}=\delta_{j}^{k}.$ There are
two interesting and useful formulas for calculating $\det[t],$ i.e.,
\begin{align}
\det[t]  &  =\underline{t}(e_{1}\wedge\ldots\wedge e_{n})\cdot(e^{1}%
\wedge\ldots\wedge e^{n}),\label{D.5a}\\
&  =\underline{t}(e^{1}\wedge\ldots\wedge e^{n})\cdot(e_{1}\wedge\ldots\wedge
e_{n}). \label{D.5b}%
\end{align}
They follow from Eq.(\ref{D.1}) by using $(e_{1}\wedge\ldots\wedge e_{n}%
)\cdot(e^{1}\wedge\ldots\wedge e^{n})=1$ which is an immediate consequence of
the formula for the Euclidean scalar product of simple $k$-vectors and the
reciprocity property of $\{e_{k}\}$ and $\{e^{k}\}$.

Each of Eq.(\ref{D.5a}) and Eq.(\ref{D.5b}) is completely equivalent to the
definition of determinant given by Eq.(\ref{D.1}).

We will end this section presenting an useful formula for the inversion of a
non-singular $(1,1)$--extensor.

Let us take $t\in ext_{1}^{1}(V)$. If $t$ is non-singular, i.e., $\det
[t]\neq0,$ then there exists its inverse $t^{-1}\in ext_{1}^{1}(V)$ which is
given by
\begin{equation}
t^{-1}(v)=\left.  \det\right.  ^{-1}[t]\underline{t}^{\dagger}(vI)I^{-1},
\label{D.6}%
\end{equation}
where $I\in\bigwedge^{n}V$ is any non-zero pseudoscalar.

To show Eq.(\ref{D.6}) we must prove that $t^{-1}$ given by the above formula
satisfies both of conditions $t^{-1}\circ t=i_{V}$ and $t\circ t^{-1}=i_{V}.$

Let $I\in\bigwedge^{n}V$ be a non-zero pseudoscalar. Take $v\in V,$ by using
the extensor identities\footnote{These extensor identities follow directly
from the fundamental identity $X\lrcorner\underline{t}(Y)=\underline
{t}(\underline{t}^{\dagger}(X)\lrcorner Y)$ with $X,Y\in\bigwedge V$. For the
first one: take $X=v,$ $Y=I$ and use $(t^{\dagger})^{\dagger}=t,$
eq.(\ref{D.1}) and $\det[t^{\dagger}]=\det[t].$ For the second one: take
$X=vI,Y=I^{-1}$ and use eq.(\ref{D.1}).} $\underline{t}^{\dagger}%
(t(v)I)I^{-1}=t(\underline{t}^{\dagger}(vI)I^{-1})=\det[t]v,$ we have that
\[
t^{-1}\circ t(v)=t^{-1}(t(v))=\left.  \det\right.  ^{-1}[t]\underline
{t}^{\dagger}(t(v)I)I^{-1}=\left.  \det\right.  ^{-1}[t]\det[t]v=i_{V}(v).
\]
And
\[
t\circ t^{-1}(v)=t(t^{-1}(v))=\left.  \det\right.  ^{-1}[t]t(\underline
{t}^{\dagger}(vI)I^{-1})=\left.  \det\right.  ^{-1}[t]\det[t]v=i_{V}(v).
\]

Finally we clarify the following. Let $T\in T^{2}(V)$ be such that for any
$u,v\in V$ we have \ $T(u,v)=t(u)\cdot v$. Then, given an arbitrary basis
$\{e_{i}\}$ of $V$ we have $T(e_{i},e_{j}):=T_{ij}=t(e_{i})\cdot e_{j}$. The
relation between the determinant of the matrix $[T_{ij}]$, denoted
\textrm{Det}$[T_{ij}]$ and $\det[t]$ is then given by:%
\begin{equation}
\mathrm{Det}[T_{ij}]=\det[t]\text{ }(e_{1}\wedge...\wedge e_{n})\cdot
(e_{1}\wedge...\wedge e_{n}). \label{relationc}%
\end{equation}
So, in general unless $\{e_{i}\}$ is an Euclidean orthonormal basis we have
that $\mathrm{Det}[T_{ij}]\neq\det[t]$.

\subsection{\ Metric and Gauge Extensors}

As we know from Section 2.8 whenever $V$ is endowed with another metric $G$
(besides $G_{E}$) there exists an unique linear operator $g$ such that the
$G$-scalar product of $X,Y\in\bigwedge V,$ namely $X\underset{g}{\cdot}Y,$ is
given by
\begin{equation}
X\underset{g}{\cdot}Y=\underline{g}(X)\cdot Y. \label{MAO.1}%
\end{equation}

Of course, \ from the definition of $(p,q)$-extensors (Section 3.1.1) we
immediately realize that the linear operator $g$ is a $(1,1)$-extensor. It is
called the pseudo-orthogonal) metric extensor for $G$.

\subsubsection{\ Metric Extensors}

We now recall that such $g\in ext_{1}^{1}(V)$ is symmetric and non-degenerate,
and has signature $(p,q),$ i.e., $g=g^{\dagger},$ $\det[g]\neq0,$ $g$ has $p$
positive and $q$ negative ($p+q=n$) eigenvalues.

Let $\{b_{j}\}$ be any orthonormal basis for $V$ with respect to $(V,G_{E}),$
i.e., $b_{j}\cdot b_{k}=\delta_{jk}.$ Once the Clifford algebra
$\mathcal{C\ell}(V,G_{E})$ has been given we are able to construct exactly $n$
(pseudo Euclidean) metric extensors with signature $(1,n-1).$ The eigenvectors
\footnote{As well-known, the eigenvalues of any orthogonal symmetric operator
are $\pm1$.} for each of them are just the basis vectors $b_{1},\ldots,b_{n}$.

Indeed, associated to $\{b_{j}\}$ we introduce the $(1,1)$-extensors
$\eta_{b_{1}},\ldots,\eta_{b_{n}}$ defined by
\begin{equation}
\eta_{b_{j}}(v)=b_{j}vb_{j}, \label{OM.1}%
\end{equation}
for each $j=1,\ldots,n$. They obviously satisfy
\begin{equation}
\eta_{b_{j}}(b_{k})=\left\{
\begin{array}
[c]{ll}%
b_{k,} & k=j\\
-b_{k}, & k\neq j
\end{array}
\right.  . \label{OM.2}%
\end{equation}
This means that $b_{j}$ is an eigenvector of $\eta_{b_{j}}$ with the
eigenvalue $+1,$ and the $n-1$ basis vectors $b_{1},\ldots,b_{j-1}%
,b_{j+1,}\ldots,b_{n}$ all are eigenvectors of $\eta_{b_{j}}$ with the same
eigenvalue $-1.$

As we can easily see, any two of these $(1,1)$-extensors commutates, i.e.,
\begin{equation}
\eta_{b_{j}}\circ\eta_{b_{k}}=\eta_{b_{k}}\circ\eta_{b_{j}},\text{ for }j\neq
k. \label{OM.3}%
\end{equation}
Moreover, they are symmetric and non-degenerate, and pseudo Euclidean
orthogonal, i.e.,
\begin{align}
\eta_{b_{j}}^{\dagger}  &  =\eta_{b_{j}}\label{OM.4a}\\
\det[\eta_{b_{j}}]  &  =(-1)^{n-1}\label{OM.4b}\\
\eta_{b_{j}}^{\dagger}  &  =\eta_{b_{j}}^{-1}. \label{OM.4c}%
\end{align}
Therefore, they all are pseudo Euclidean orthogonal metric extensors with
signature $(1,n-1)$.

The extended of $\eta_{b_{j}}$ is given by
\begin{equation}
\underline{\eta_{b_{j}}}(X)=b_{j}Xb_{j}, \label{OM.5}%
\end{equation}

\subsubsection{Constructing Metrics of Signature $(p,n-p)$}

We can now construct a pseudo orthogonal metric operator with signature
$(p,n-p)$ and whose orthonormal eigenvectors are just the basis vectors
$b_{1},\ldots,b_{n}$. It is defined by
\begin{equation}
\eta_{b}=(-1)^{p+1}\eta_{b_{1}}\circ\cdots\circ\eta_{b_{p}}, \label{OM.6}%
\end{equation}
i.e.,
\begin{equation}
\eta_{b}(a)=(-1)^{p+1}b_{1}\ldots b_{p}ab_{p}\ldots b_{1}. \label{OM.7}%
\end{equation}

Wee have that
\begin{equation}
\eta_{b}(b_{k})=\left\{
\begin{array}
[c]{ll}%
b_{k}, & k=1,\ldots,p\\
-b_{k}, & k=p+1,\ldots,n
\end{array}
\right.  , \label{OM.8}%
\end{equation}
which means that $b_{1},\ldots,b_{p}$ are eigenvectors of $\eta_{b}$ with the
same eigenvalue $+1$, and $b_{p+1},\ldots,b_{n}$ are eigenvectors of $\eta
_{b}$ with the same eigenvalue $-1.$

The extensor $\eta_{b}$ is symmetric and non-degenerate, and orthogonal, i.e.,
$\eta_{b}^{\dagger}=\eta_{b}$, $\det[\eta_{b}]=(-1)^{n-p}$ and $\eta
_{b}^{\dagger}=\eta_{b}^{-1}$ and thus $\eta_{b}$ is a pseudo orthogonal
metric extensor with signature $(p,n-p)$.

The extended of $\eta_{b}$ is obviously given by
\begin{equation}
\underline{\eta_{b}}(X)=(-1)^{p+1}b_{1}\ldots b_{p}Xb_{p}\ldots b_{1}.
\label{OM.10}%
\end{equation}

What is the most general pseudo-orthogonal metric extensor with signature
$(p,n-p)$?

To find the answer, let $\eta$ be any pseudo-orthogonal metric extensor with
signature $(p,n-p).$ The symmetry of $\eta$ implies the existence of exactly
$n$ Euclidean orthonormal eigenvectors $u_{1},\ldots,u_{n}$ for $\eta$ which
form just a basis for $V$. Since $\eta$ is pseudo-orthogonal and its signature
is $(p,n-p)$, it follows that the eigenvalues of $\eta$ are equal $\pm1$ and
the eigenvalues equation for $\eta$ can be written (re-ordering $u_{1}%
,\ldots,u_{n}$ if necessary) as
\[
\eta(u_{k})=\left\{
\begin{array}
[c]{ll}%
u_{k}, & k=1,\ldots,p\\
-u_{k}, & k=p+1,\ldots,n
\end{array}
\right.  .
\]

Now, due to the orthonormality of both $\{b_{k}\}$ and $\{u_{k}\},$ there must
be an orthogonal $(1,1)$-extensor $\Theta$ such that $\Theta(b_{k})=u_{k},$
for each $k=1,\ldots,n,$ i.e., for all $a\in V:\Theta(a)=\overset{n}%
{\underset{j=1}{%
%TCIMACRO{\dsum }%
%BeginExpansion
{\displaystyle\sum}
%EndExpansion
}}(a\cdot b_{j})u_{j}.$

Then, we can write
\begin{align*}
\Theta\circ\eta_{b}\circ\Theta^{\dagger}(u_{k})  &  =\Theta\circ\eta_{b}%
(b_{k})=\Theta(\left\{
\begin{array}
[c]{ll}%
b_{k}, & k=1,\ldots,p\\
-b_{k}, & k=p+1,\ldots,n
\end{array}
\right.  )\\
&  =\left\{
\begin{array}
[c]{ll}%
u_{k}, & k=1,\ldots,p\\
-u_{k}, & k=p+1,\ldots,n
\end{array}
\right.  =\eta(u_{k}),
\end{align*}
for each $k=1,\ldots,n.$ Thus, we have
\begin{equation}
\eta=\Theta\circ\eta_{b}\circ\Theta^{\dagger}. \label{OM.11}%
\end{equation}
By putting Eq.(\ref{OM.6}) into Eq.(\ref{OM.11}) we get
\begin{equation}
\eta=(-1)^{p+1}\eta_{1}\circ\cdots\circ\eta_{p}, \label{OM.12}%
\end{equation}
where each of $\eta_{j}\equiv\Theta\circ\eta_{b_{j}}\circ\Theta^{\dagger}$ is
an Euclidean orthogonal metric extensor with signature $(1,n-p)$.

But, by using the vector identity $abc=(a\cdot b)c-(a\cdot c)b+(b\cdot
c)a+a\wedge b\wedge c,$ with $a,b,c\in V,$ we can prove that
\begin{equation}
\eta_{j}(v)=\Theta(b_{j})v\Theta(b_{j}), \label{OM.13}%
\end{equation}
for each $j=1,\ldots,p$.

Now, by using Eq.(\ref{OM.13}) we can write Eq.(\ref{OM.12}) in the remarkable
form
\begin{equation}
\eta(v)=(-1)^{p+1}\underline{\Theta}(b_{1}\ldots b_{p})v\underline{\Theta
}(b_{p}\ldots b_{1}). \label{OM.14}%
\end{equation}

\subsubsection{Metric Adjoint Extensor}

Given a metric extensor $g$ for $V$, to each $t\in1$-$ext(\bigwedge
_{1}^{\diamond}V;\bigwedge_{2}^{\diamond}V)$ we can assign $t^{\dagger(g)}%
\in1$-$ext(\bigwedge_{2}^{\diamond}V;\bigwedge_{1}^{\diamond}V)$ defined as
follows:
\begin{equation}
t^{\dagger(g)}=\underline{g}^{-1}\circ t^{\dagger}\circ\underline{g}.
\label{MAU.2}%
\end{equation}
It will be called the\emph{\ metric adjoint} of $t.$

As we can easily see, $t^{\dagger(g)}$ is the unique extensor from
$\bigwedge_{2}^{\diamond}V$ to $\bigwedge_{1}^{\diamond}V$ which satisfies the
fundamental property
\begin{equation}
X\underset{g}{\cdot}t^{\dagger(g)}(Y)=t(X)\underset{g}{\cdot}Y, \label{MAO.3}%
\end{equation}
for all $X\in\bigwedge_{1}^{\diamond}V$ and $Y\in\bigwedge_{2}^{\diamond}V$
and where we recall the symbol $\underset{g}{\cdot}$ reffers to the scalar
product determined by a general (nondegenrated) metric $G$ or equivalently its
associated extensor $g$ according to Eq.(\ref{MAO.1}).

The noticeable property given by Eq.(\ref{MAO.3}) is the \emph{metric version}
of the fundamental property given by Eq.(\ref{SAO.3}).

\paragraph{Lorentz Extensor}

A $(1,1)$-extensor over $V,$ namely $\Lambda,$ is said to be $\eta$-orthogonal
if and only if for all $v,w\in V$
\begin{equation}
\Lambda(v)\underset{\eta}{\cdot}\Lambda(w)=v\underset{\eta}{\cdot}w.
\label{OE.1}%
\end{equation}

Recalling the non-degeneracy of the $\eta$-scalar product, Eq.(\ref{OE.1})
implies that
\begin{equation}
\Lambda^{\dagger(\eta)}=\Lambda^{-1}. \label{OE.2}%
\end{equation}
Or, by taking into account Eq.(\ref{MAU.2}), we can still write
\begin{equation}
\Lambda^{\dagger}\circ\eta\circ\Lambda=\eta. \label{OE.3}%
\end{equation}

We emphasize that the $\eta$-scalar product condition given by Eq.(\ref{OE.1})
is logically equivalent to each of Eq.(\ref{OE.2}) and Eq.(\ref{OE.3}).

Sometimes, when the signature of the $\eta$-orthogonal $(1,1)$-extensor is
$(1,n-1)$, the extensor $\Lambda$ is called a \emph{Lorentz extensor}.

\subsubsection{Gauge Extensors}

\label{ONE}Let $g$ and $\eta$ be pseudo-orthogonal metric extensors, of the
same signature $(p,n-p).$ Then, we can show that there exists a non-singular
$(1,1)$-extensor $h$ such that
\begin{equation}
g=h^{\dagger}\circ\eta\circ h. \label{GE.1}%
\end{equation}

Such $h$ is given by
\begin{equation}
h=d_{\sigma}\circ d_{\sqrt{\left\vert \lambda\right\vert }}\circ\Theta_{uv},
\label{GE.3}%
\end{equation}
where $d_{\sigma}$ is a pseudo-orthogonal metric extensor, $d_{\sqrt
{\left\vert \lambda\right\vert }}$ is a metric extensor, and $\Theta_{uv}$ is
a \ pseudo-orthogonal operator which are defined by
\begin{align}
d_{\sigma}(a)  &  =\overset{n}{\underset{j=1}{\sum}}\sigma_{j}(a\cdot
u_{j})u_{j}\label{GE.3a}\\
d_{\sqrt{\left\vert \lambda\right\vert }}(a)  &  =\overset{n}{\underset
{j=1}{\sum}}\sqrt{\left\vert \lambda_{j}\right\vert }(a\cdot u_{j}%
)u_{j}\label{GE.3b}\\
\Theta_{uv}(a)  &  =\overset{n}{\underset{j=1}{\sum}}(a\cdot u_{j})v_{j},
\label{GE.3c}%
\end{align}
where $\sigma_{1},\ldots,\sigma_{n}$ are real numbers with $\sigma_{1}%
^{2}=\cdots=\sigma_{n}^{2}=1,$ $\lambda_{1},\ldots,\lambda_{n}$ are the
eigenvalues of $g,$ and $u_{1},\ldots,u_{n}$ and $v_{1},\ldots,v_{n}$ are
respectively the orthonormal eigenvectors of $\eta$ and $g$. For a proof see,
e.g., \cite{rodoliv2006}.

\subsubsection{ Golden Formula}

Let $h$ be any gauge extensor for $g,$ i.e., $g=h^{\dagger}\circ\eta\circ h,$
where $\eta$ is a pseudo-orthogonal metric extensor with the same signature as
$g.$ Let $\underset{g}{\ast}$ mean either $\wedge$ (exterior product),
$\underset{g}{\cdot}$ ($g$-scalar product), $\underset{g}{\lrcorner}%
,\underset{g}{\llcorner}$ ($g$-contracted products) or $\underset{g}{}$
($g$-Clifford product). And analogously for $\underset{\eta}{\ast}$.

The $g$-metric products $\underset{g}{\ast}$ and the $\eta$-metric products
are related by a remarkable formula, called in what follows the \emph{golden
formula}. For all $X,Y\in\bigwedge V$
\begin{equation}
\underline{h}(X\underset{g}{\ast}Y)=[\underline{h}(X)\underset{\eta}{\ast
}\underline{h}(Y)], \label{GF.1}%
\end{equation}
where $\underline{h}$ denotes the \textit{extended }\cite{10} of $h$. A proof
of the golden formula can be found in, e.g., \cite{rodoliv2006}.

\subsection{Hodge Extensors}

\subsubsection{Standard Hodge Extensor}

Let $\{e_{j}\}$ and $\{e^{j}\}$ be two Euclidean reciprocal bases to each
other for $V,$ i.e., $e_{j}\cdot e^{k}=\delta_{j}^{k}.$ Associated to them we
define a non-zero pseudoscalar
\begin{equation}
\tau=\sqrt{e_{\wedge}\cdot e_{\wedge}}e^{\wedge}, \label{SHE.1}%
\end{equation}
where $e_{\wedge}\equiv e_{1}\wedge\ldots\wedge e_{n}\in\bigwedge^{n}V$ and
$e^{\wedge}\equiv e^{1}\wedge\ldots\wedge e^{n}\in\bigwedge^{n}V$. Note that
$e_{\wedge}\cdot e_{\wedge}>0,$ since an Euclidean scalar product is positive
definite. Such $\tau$ will be called a \emph{standard volume pseudoscalar} for
$V.$

The standard volume pseudoscalar has the fundamental property
\begin{equation}
\tau\cdot\tau=\tau\lrcorner\widetilde{\tau}=\tau\widetilde{\tau}=1,
\label{SHE.2}%
\end{equation}
which follows from the obvious result $e_{\wedge}\cdot e^{\wedge}=1$.

From Eq.(\ref{SHE.2}), we can easily get an expansion formula for
pseudoscalars of $\bigwedge^{n}V,$ i.e.,
\begin{equation}
I=(I\cdot\tau)\tau. \label{SHE.3}%
\end{equation}

The extensor $\star\in ext(V)$ which is defined by $\star:\bigwedge
V\rightarrow\bigwedge V$ such that
\begin{equation}
\star X=\widetilde{X}\lrcorner\tau=\widetilde{X}\tau, \label{SHE.4}%
\end{equation}
will be called a \emph{standard Hodge extensor }on $V$.

It should be noticed that
\begin{equation}
\text{if }X\in\bigwedge^{p}V,\text{ then }\star X\in\bigwedge^{n-p}V.
\label{SHE.4a}%
\end{equation}
That means that $\star$ can be also defined as a $(p,n-p)$-extensor over $V. $

The extensor over $V,$ namely $\star^{-1},$ which is given by $\star
^{-1}:\bigwedge V\rightarrow\bigwedge V$ such that
\begin{equation}
\star^{-1}X=\tau\llcorner\widetilde{X}=\tau\widetilde{X} \label{SHE.5}%
\end{equation}
is the \emph{inverse extensor} of $\star$. We have immediately that for any
$X\in\bigwedge V$ $\star^{-1}\circ\star X=\tau\widetilde{\tau}X=X$,
$\star\circ\star^{-1}X=X\widetilde{\tau}\tau=X$, i.e., $\star^{-1}\circ
\star=\star\circ\star^{-1}=i_{\bigwedge V},$ where $i_{\bigwedge V}\in ext(V)$
is the so-called \emph{identity function }for $\bigwedge V.$

We recall now some identities that we shall use in the next papers of this series.

(i) Let for $X,Y\in\bigwedge V$. We have
\begin{equation}
(\star X)\cdot(\star Y)=X\cdot Y, \label{SHE.6}%
\end{equation}
which means that the standard Hodge extensor preserves the Euclidean scalar product.

(ii) Let us take $X,Y\in\bigwedge^{p}V.$ By using Eq.(\ref{SHE.3}) together
with the multivector identity $(X\wedge Y)\cdot Z=Y\cdot(\widetilde
{X}\lrcorner Z),$ and Eq.(\ref{SHE.6}) we get
\begin{equation}
X\wedge(\star Y)=(X\cdot Y)\tau. \label{SHE.7}%
\end{equation}
This noticeable identity is completely equivalent to the definition of the
standard Hodge extensor given by Eq.(\ref{SHE.4}) and indeed it is analogous
to the one used to define the Hodge dual of form fields in texts dealing with
the Cartan calculus of differential forms.

(iii) Let $X\in\bigwedge^{p}V$ and $Y\in\bigwedge^{n-p}V.$ By using the
multivector identity $(X\lrcorner Y)\cdot Z=Y\cdot(\widetilde{X}\wedge Z)$ and
Eq.(\ref{SHE.3}) we get
\begin{equation}
(\star X)\cdot Y\tau=X\wedge Y. \label{SHE.8}%
\end{equation}

\subsubsection{Metric Hodge Extensor}

Our objective in this section is to find \ a formula which permit us to
related the Hodge dual operators associated to two distinct metrics. Let then,
$g$ be a metric extensor on $V$ with signature $(p,q),$ i.e., $g\in
ext_{1}^{1}(V)$ such that $g=g^{\dagger}$ and $\det[g]\neq0,$ and it has $p$
positive and $q$ negative eigenvalues. Associated to $\{e_{j}\}$ and
$\{e^{j}\}$ we can define another non-zero pseudoscalar
\begin{equation}
\underset{g}{\tau}\text{ }=\sqrt{\left\vert e_{\wedge}\underset{g}{\cdot
}e^{\wedge}\right\vert }e^{\wedge}=\sqrt{\left\vert \det[g]\right\vert }\tau.
\label{MHE.1}%
\end{equation}
It will be called a metric volume pseudoscalar for $V.$ It has the fundamental
property
\begin{equation}
\underset{g}{\tau}\underset{g^{-1}}{\cdot}\underset{g}{\tau}=\underset{g}%
{\tau}\underset{g^{-1}}{\lrcorner}\widetilde{\underset{g}{\tau}}=\underset
{g}{\tau}\underset{g^{-1}}{}\widetilde{\underset{g}{\tau}}=(-1)^{q}.
\label{MHE.2}%
\end{equation}
Eq.(\ref{MHE.2}) follows from Eq.(\ref{SHE.2}) by taking into account the
definition of determinant of a $(1,1)$-extensor, and recalling that
$sgn(\det[g])=(-1)^{q}.$

An expansion formula for pseudoscalars of $\bigwedge^{n}V$ can be also
obtained from Eq.(\ref{MHE.2}), i.e.,
\begin{equation}
I=(-1)^{q}(I\underset{g^{-1}}{\cdot}\underset{g}{\tau})\underset{g}{\tau}.
\label{MHE.3}%
\end{equation}

The extensor $\underset{g}{\star}\in ext(V)$ which is defined by $\underset
{g}{\star}:\bigwedge V\rightarrow\bigwedge V$ such that
\begin{equation}
\underset{g}{\star}X=\widetilde{X}\underset{g^{-1}}{\lrcorner}\underset
{g}{\tau}=\widetilde{X}\underset{g^{-1}}{}\underset{g}{\tau} \label{MHE.4}%
\end{equation}
will be called a \emph{metric Hodge extensor }on $V.$ It should be noticed
that in general we need to use of both the $g$ and $g^{-1}$ metric Clifford algebras.

We see that
\begin{equation}
\text{if }X\in\bigwedge^{p}V,\text{ then }\underset{g}{\star}X\in
\bigwedge^{n-p}V. \label{MHE.4a}%
\end{equation}
It means that $\underset{g}{\star}\in ext(V)$ can also be defined as
$\underset{g}{\star}\in ext_{p}^{n-p}(V).$

The extensor over $V,$ namely $\underset{g}{\star}^{-1},$ which is given by
$\underset{g}{\star}^{-1}:\bigwedge V\rightarrow\bigwedge V$ such that
\begin{equation}
\underset{g}{\star}^{-1}X=(-1)^{q}\underset{g}{\tau}\underset{g^{-1}%
}{\llcorner}\widetilde{X}=(-1)^{q}\underset{g}{\tau}\underset{g^{-1}}%
{}\widetilde{X} \label{MHE.5}%
\end{equation}
is the \emph{inverse extensor} of $\underset{g}{\star}.$

Let us take $X\in\bigwedge V.$ By using Eq.(\ref{MHE.2}), we have indeed that
$\underset{g}{\star}^{-1}\circ\underset{g}{\star}X=(-1)^{q}\underset{g}{\tau
}\underset{g^{-1}}{}\widetilde{\underset{g}{\tau}}\underset{g^{-1}}{}X=X,$ and
$\underset{g}{\star}\circ\underset{g}{\star}^{-1}X=(-1)^{q}X\underset{g^{-1}%
}{}\widetilde{\underset{g}{\tau}}\underset{g^{-1}}{}\underset{g}{\tau}=X,$
i.e., $\underset{g}{\star}^{-1}\circ\star=\underset{g}{\star}\circ\underset
{g}{\star}^{-1}=i_{\bigwedge V}.$

Take $X,Y\in\bigwedge V.$ The identity $(X\underset{g^{-1}}{}A)\underset
{g^{-1}}{\cdot}Y=X\underset{g^{-1}}{\cdot}(Y\underset{g^{-1}}{}\widetilde{A})$
and Eq.(\ref{MHE.2}) yield
\begin{equation}
(\underset{g}{\star}X)\underset{g^{-1}}{\cdot}(\underset{g}{\star}%
Y)=(-1)^{q}X\underset{g^{-1}}{\cdot}Y. \label{MHE.6}%
\end{equation}

Take $X,Y\in\bigwedge^{p}V.$ Eq.(\ref{MHE.3}), the identity $(X\wedge
Y)\underset{g^{-1}}{\cdot}Z=Y\underset{g^{-1}}{\cdot}(\widetilde{X}%
\underset{g^{-1}}{\lrcorner}Z)$ and Eq.(\ref{MHE.6}) allow us to obtain
\begin{equation}
X\wedge(\underset{g}{\star}Y)=(X\underset{g^{-1}}{\cdot}Y)\underset{g}{\tau}.
\label{MHE.7}%
\end{equation}
This remarkable property is completely equivalent to the definition of the
metric Hodge extensor given by Eq.(\ref{MHE.4}).

Take $X\in\bigwedge^{p}V$ and $Y\in\bigwedge^{n-p}V.$ The use of identity
$(X\underset{g^{-1}}{\lrcorner}Y)\underset{g^{-1}}{\cdot}Z=Y\underset{g^{-1}%
}{\cdot}(\widetilde{X}\wedge Z)$ and Eq.(\ref{MHE.3}) yield
\begin{equation}
(\underset{g}{\star}X)\underset{g^{-1}}{\cdot}Y\underset{g}{\tau}%
=(-1)^{q}X\wedge Y. \label{MHE.8}%
\end{equation}

It might as well be asked what is the relationship between the standard and
metric Hodge extensors as defined above by Eq.(\ref{SHE.4}) and
Eq.(\ref{MHE.4}).

Take $X\in\bigwedge V.$ By using Eq.(\ref{MHE.1}), the multivector identity
for an invertible $(1,1)$-extensor $\underline{t}^{-1}(X)\lrcorner
Y=\underline{t}^{\dagger}(X\lrcorner\underline{t}^{\ast}(Y)),$ and the
definition of determinant of a $(1,1)$-extensor we have that
\begin{align*}
\underset{g}{\star}X  &  =\underline{g}^{-1}(\widetilde{X})\lrcorner
\sqrt{\left\vert \det[g]\right\vert }\tau=\sqrt{\left\vert \det[g]\right\vert
}\underline{g}(\widetilde{X}\lrcorner\underline{g}^{-1}(\tau))\\
&  =\frac{\sqrt{\left\vert \det[g]\right\vert }}{\det[g]}\underline
{g}(\widetilde{X}\lrcorner\tau)=\frac{sgn(\det[g])}{\sqrt{\left\vert
\det[g]\right\vert }}\underline{g}\circ\star(X),
\end{align*}
i.e.,
\begin{equation}
\underset{g}{\star}=\frac{(-1)^{q}}{\sqrt{\left\vert \det[g]\right\vert }%
}\underline{g}\circ\star. \label{MHE.9}%
\end{equation}
Eq.(\ref{MHE.9}) is then the formula which relates the metric Hodge extensor
$\underset{g}{\star}$ with the standard Hodge extensor $\star.$

We already know that for any metric extensor $g\in ext_{1}^{1}(V)$ there
exists a non-singular $(1,1)$-extensor $h\in ext_{1}^{1}(V)$ (the \emph{gauge
extensor} for $g$) such that
\begin{equation}
g=h^{\dagger}\circ\eta\circ h, \label{MHE.10}%
\end{equation}
where $\eta\in ext_{1}^{1}(V)$ is a pseudo-orthogonal metric extensor with the
same signature as $g$.

We now obtain a noticiable formula which relates the $g$-metric Hodge extensor
with the $\eta$-metric Hodge extensor.

As we know, the $g$ and $g^{-1}$ contracted products $\underset{g}{\lrcorner}$
and $\underset{g^{-1}}{\lrcorner}$ are related to the $\eta$-contracted
product$\underset{\eta}{\lrcorner}$ (recall that $\eta=\eta^{-1}$) by the
following \emph{golden formulas}
\begin{align}
\underline{h}(X\underset{g}{\lrcorner}Y)  &  =\underline{h}(X)\underset{\eta
}{\lrcorner}\underline{h}(Y),\label{MHE.10a}\\
\underline{h}^{*}(X\underset{g^{-1}}{\lrcorner}Y)  &  =\underline{h}%
^{*}(X)\underset{\eta}{\lrcorner}\underline{h}^{*}(Y). \label{MHE.10b}%
\end{align}

Now, take $X\in\bigwedge V.$ By using Eq.(\ref{MHE.10b}), Eq.(\ref{MHE.1}),
the definition of determinant of a $(1,1)$-extensor, Eq.(\ref{MHE.10}) and the
obvious equation $\underset{\eta}{\tau}=\tau$ we have that
\begin{align*}
\underset{g}{\star}X  &  =\underline{h}^{\dagger}(\underline{h}^{\ast
}(\widetilde{X})\underset{\eta}{\lrcorner}\underline{h}^{\ast}(\underset
{g}{\tau}))=\sqrt{\left\vert \det[g]\right\vert }\underline{h}^{\dagger
}(\underline{h}^{\ast}(\widetilde{X})\underset{\eta}{\lrcorner}\det[h^{\ast
}]\tau)\\
&  =\left\vert \det[h]\right\vert \det[h^{\ast}]\underline{h}^{\dagger
}(\widetilde{\underline{h}^{\ast}(X)}\underset{\eta^{-1}}{\lrcorner}%
\underset{\eta}{\tau})=sgn(\det[h])\underline{h}^{\dagger}\circ\underset{\eta
}{\star}\circ\underline{h}^{\ast}(X),
\end{align*}
i.e.,
\begin{equation}
\underset{g}{\star}\text{ }=sgn(\det[h])\underline{h}^{\dagger}\circ
\underset{\eta}{\star}\circ\underline{h}^{\ast}. \label{MHE.11}%
\end{equation}
This formula which relates the $g$-metric Hodge extensor $\underset{g}{\star}$
with the $\eta$-metric Hodge extensor $\underset{\eta}{\star}$ will play an
important role in the applications we have in mind.

\subsection{Some Useful Bases}

For future reference we end this paper introducing some related bases that are
useful in the developments that we have in mind in this series.

(i) Let $\{e_{k}\}$ be any basis for $V,$ and $\{e^{k}\}$ be its Euclidean
reciprocal basis for $V$, i.e., $e_{k}\cdot e^{l}=\delta_{k}^{l}.$ Let us take
a non-singular $(1,1)$-extensor $\lambda.$ Then, it is easily seen that the
$n$ vectors $\lambda(e_{1}),\ldots,\lambda(e_{n})\in V$ and the $n$
vectors\footnote{Recall that $\lambda^{*}=(\lambda^{-1})^{\dagger}%
=(\lambda^{\dagger})^{-1}.$} $\lambda^{\ast}(e^{1}),\ldots,\lambda^{\ast
}(e^{n})\in V$ define two well-defined Euclidean reciprocal bases for $V,$
i.e.,
\begin{equation}
\lambda(e_{k})\cdot\lambda^{\ast}(e^{l})=\delta_{k}^{l}. \label{RT.1}%
\end{equation}
The bases $\{\lambda(e_{k})\}$ and $\{\lambda^{\ast}(e^{k})\}$ are
conveniently said to be a $\lambda$\emph{-deformation} of the bases
$\{e_{k}\}$ and $\{e^{k}\}.$ Sometimes, the first ones are named as the
$\lambda$\emph{-deformed bases }of the second ones.

(ii) Let $h$ be a gauge extensor for $g,$ and $\eta$ be a pseudo-orthogonal
metric extensor with the same signature as $g.$ According to Eq.(\ref{GE.1})
the $g$-scalar product and $g^{-1}$-scalar product are related to the $\eta
$-scalar product by the following formulas
\begin{equation}
X\underset{g}{\cdot}Y=\underline{h}(X)\underset{\eta}{\cdot}\underline
{h}(Y)\text{, }X\underset{g^{-1}}{\cdot}Y=\underline{h}^{\ast}(X)\underset
{\eta}{\cdot}\underline{h}^{\ast}(Y). \label{GB.1b}%
\end{equation}

The $\eta$-deformed bases $\{h(e_{k})\}$ and $\{h^{\ast}(e^{k})\}$ satisfy the
noticeable properties\footnote{Recall that $g_{jk}\equiv g(e_{j})\cdot
e_{k}=G(e_{j},e_{k})\equiv G_{jk}$ and $g^{jk}\equiv g^{-1}(e^{j})\cdot
e^{k}=G^{jk}$ are the $jk$-entries of the inverse matrix for $[G_{jk}]$.}
\begin{equation}
h(e_{j})\underset{\eta}{\cdot}h(e_{k})=g\text{, }h^{\ast}(e^{j})\underset
{\eta}{\cdot}h^{\ast}(e^{k})=g^{jk}. \label{GB.2b}%
\end{equation}
The bases $\{h(e_{k})\}$ and $\{h^{\ast}(e^{k})\}$ are called the \emph{gauge
bases} associated with $\{e_{k}\}$ and $\{e^{k}\}.$

(iii) Let $u_{1},\ldots,u_{n}$ be the $n$ \emph{Euclidean} orthonormal
eigenvectors of $\eta$, i.e., the eigenvalues equation for $\eta$ can be
written (reordering $u_{1},\ldots,u_{n}$ if necessary) as
\[
\eta(u_{k})=\left\{
\begin{array}
[c]{cc}%
u_{k}, & k=1,\ldots,p\\
-u_{k}, & k=p+1,\ldots,n
\end{array}
\right.  ,
\]
and $u_{j}\cdot u_{k}=\delta_{jk}.$

The $h^{-1}$-deformed bases $\{h^{-1}(u_{k})\}$ and $\{h^{\dagger}(u_{k})\}$
satisfy the remarkable properties
\begin{equation}
h^{-1}(u_{j})\underset{g}{\cdot}h^{-1}(u_{k})=\eta_{jk}\text{, }h^{\dagger
}(u_{j})\underset{g^{-1}}{\cdot}h^{\dagger}(u_{k})=\eta_{jk}, \label{TB.1b}%
\end{equation}
where

$\eta_{jk}\equiv\eta(u_{j})\cdot u_{k}=\left\{
\begin{array}
[c]{cc}%
1, & j=k=1,\ldots,p\\
-1, & j=k=p+1,\ldots,n\\
0, & j\equiv k
\end{array}
\right.  .$

The bases $\{h^{-1}(u_{k})\}$ and $\{h^{\dagger}(u_{k})\}$ are called the
\emph{tetrad bases} associated with $\{u_{k}\}$.

\section{Conclusions}

This paper, the first in a series of four presents the theory of geometric
(Clifford) algebras and the theory of extensors. Our presentation has been
devised in order to provide a powerful computational tool which permits its
efficient application in the study of differential geometry (in the next
papers of the series) in a natural and simple way. Besides a thoughtful
presentation of the concepts, we detailed many calculations in order to
provide the \textquotedblleft tricks of the trade\textquotedblright\ to
readers interested in applications\footnote{We advise that someone interested
in aspects of the general theory of Clifford algebras not covered here we
strongly recommend the excellent textbooks \cite{1,12,13,14}.}. Among the
results obtained worth to quote here, a distinction goes to the theory of
deformation of the Euclidean Clifford algebra $\mathcal{C}\ell(V,G_{E})$
through the use of metric and gauge extensors which permits to generate all
other Clifford algebras $\mathcal{C}\ell(V,G)$ (where $G$ is a metric of
signature $(p,q)$ with $p+q=n=\dim V$ and $G_{E}$ is an Euclidean metric in
$V$) and derivation of the golden formula, essential for reducing all
calculations in any Clifford algebra $\mathcal{C}\ell(V,G)$ to the ones in
$\mathcal{C}\ell(V,G_{E})$, thus providing many useful formulas as, e.g., a
remarkable relation between Hodge (star) operators associated to $G$ and
$G_{E}$.

\textbf{Acknowledgments: } V. V. Fern\'{a}ndez and A. M. Moya are very
grateful to Mrs. Rosa I. Fern\'{a}ndez who gave to them material and spiritual
support at the starting time of their research work. This paper could not have
been written without her inestimable help. Authors are also grateful to Drs.
E. Notte-Cuello and E. Capelas de Oliveira for valuable discussions.


\begin{thebibliography}{99}                                                                                               %


\bibitem {1}{\footnotesize Ablamowicz, R., Baylis, W. E., Branson, T.,
Lounesto, P., Porteous, I., Ryan, J., Selig, J.M., Sobczyk, G., in Ablamowicz,
R. and Sobczyk, G. (eds.), \textit{Lectures on Clifford (Geometrical) Algebras
and Applications}, Birkh\"{a}user, Boston (2004).}

\bibitem {3}{\footnotesize Moya, A. M., Fern\'{a}ndez, V. V., and Rodrigues,
W. A., Jr. Multivector and Extensor Fields on Arbitrary Manifolds,
\textit{Int. J. Geom. Meth. Mod. Phys. }\textbf{4 }(6) (2007).}

\bibitem {10}{\footnotesize Fern\'{a}ndez, V. V., Moya, A. M., and Rodrigues,
W. A., Jr., Mutivector and Extensor Calculus, Special Issue of \textit{Adv. in
Appl. Clifford Algebras} \textbf{11}(S3), 1-103\ (2001)}.

\bibitem {8}{\footnotesize Hestenes, D., and Sobcyk, G., \textit{Clifford
Algebras to Geometrical Calculus}, D. Reidel Publ. Co., Dordrecht, 1984.}

\bibitem {9}{\footnotesize Lasenby, A., Doran, C., and Gull, S., Gravity,
Gauge Theories and Geometric Algebras, \textit{Phil. Trans. R. Soc.}
\textbf{356}, 487-582 (1998).}

\bibitem {12}{\footnotesize Lounesto, P., \emph{Clifford Algebras and
Spinors}, London Math. Soc., Lecture Notes Series 239, Cambridge University
Press, Cambridge, 1997; second edition 2001.}

\bibitem {moro}{\footnotesize Mosna, R. A. and Rodrigues, W. A. Jr., The
Bundles of Algebraic and Dirac-Hestenes Spinor Fields, \textit{J. Math. Phys.}
\textbf{45}, 2945-2966 (2004).}

\bibitem {4}{\footnotesize Fern\'{a}ndez, V. V., Moya, A. M., and Rodrigues,
W. A., Jr., Geometric and Extensor Algebras and the Differential Geometry of
Arbitrary Manifolds, \textit{Int. J. Geom. Meth. Mod. Phys. }\textbf{4 }(7)
(2007).}

\bibitem {7}{\footnotesize Fern\'{a}ndez, V. V., Moya, A. M., da Rocha, R.,
and Rodrigues, W. A., Jr., Clifford and Extensor Calculus and the Riemann and
Ricci Fields of Deformed Structures }$(M,\nabla,\mathbf{\eta})$%
{\footnotesize and }$(M,\nabla,%
%TCIMACRO{\TeXButton{slg}{\slg}}%
%BeginExpansion
\slg
%EndExpansion
)${\footnotesize , \textit{Int. J. Geom. Meth. Mod. Phys. }\textbf{4} (7)}
{\footnotesize \textit{ }(2007).}

\bibitem {13}{\footnotesize Porteous, I. R., \textit{Clifford Algebras and the
Classical Groups}, Cambridge Studies in Advanced Mathematics 50, Cambridge
University Press, Cambridge, 1995; second edition 2001.}

\bibitem {14}{\footnotesize Porteous, I. R., \textit{Topological Geometry},
Van Nostrand Reinhold, London, 1969; 2nd edition, Cambridge University Press,
Cambridge, 1981.}

\bibitem {rod041}{\footnotesize Rodrigues, W. A. Jr., Algebraic and
Dirac-Hestenes Spinors and Spinor Fields, \textit{J. Math. Phys. }\textbf{45},
2908-2944 (2004).}

\bibitem {rodoliv2006}{\footnotesize Rodrigues, W. A. Jr., and Oliveira, E.
Capelas, \textit{The Many Faces of Maxwell, Dirac and Einstein Equations. A
Clifford Bundle Approach}, Lecture Notes in Physics \textbf{722}, Springer,
New York, 2007.}

\bibitem {11}{\footnotesize Sobczyk, G., Direct Integration, in Baylis, W. E.
(ed.), \textit{Clifford (Geometric) Algebras with Applications in Physics,
Mathematics and Engineering}, pp. 53-64, Birkh\"{a}user, Boston, 1999.}
\end{thebibliography}
\end{document}